\documentclass[11pt,a4paper]{article}
\usepackage{latexsym}
\usepackage{amsmath}
\usepackage{amssymb}
\usepackage{amscd}
\usepackage{amsthm}
\usepackage[all]{xy}

\newtheorem{thm}{Theorem}[section]
\newtheorem{cor}[thm]{Corollary}
\newtheorem{lem}[thm]{Lemma}
\newtheorem{prop}[thm]{Proposition}

\theoremstyle{definition}
\newtheorem{defn}[thm]{Definition}
\newtheorem{rem}[thm]{Remark}
\newtheorem{exa}[thm]{Example}

\DeclareMathOperator{\GL}{\mathbf{GL}}

\newcommand{\R}{\mathbb R}
\newcommand{\Z}{\mathbb Z}

\newcommand{\C}{\mathbb C}

\newif\ifpdf \pdftrue
\ifx\pdfoutput\undefined \pdffalse \fi \ifx\pdfoutput\relax
\pdffalse \fi

\usepackage{makeidx}
\makeindex

\begin{document}

\title{Uniqueness of equivariant singular Bott-Chern classes}

\author{Shun Tang}

\date{}

\maketitle

\vspace{-10mm}

\hspace{5cm}\hrulefill\hspace{5.5cm} \vspace{5mm}

\textbf{Abstract.} In this paper, we shall discuss possible
theories of defining equivariant singular Bott-Chern classes and
corresponding uniqueness property. By adding a natural axiomatic
characterization to the usual ones of equivariant Bott-Chern
secondary characteristic classes, we will see that the
construction of Bismut's equivariant Bott-Chern singular currents
provides a unique way to define a theory of equivariant singular
Bott-Chern classes. This generalizes J. I. Burgos Gil and R.
Li\c{t}canu's discussion to the equivariant case. As a byproduct
of this study, we shall prove a concentration formula which can be
used to prove an arithmetic concentration theorem in Arakelov
geometry.

\textbf{R\'{e}sum\'{e}.} Dans cet article, nous allons discuter
les th\'{e}ories \'{e}ventuelles de d\'{e}finir les classes de
Bott-Chern \'{e}quivariantes singuli\`{e}res et la
propri\'{e}t\'{e} d'unicit\'{e} correspondante. En ajoutant une
caract\'{e}risation axiomatique naturelle \`{a} lesquelles
habituelles des classes caract\'{e}ristiques secondaires de
Bott-Chern \'{e}quivariantes, nous verrons que la construction des
courants de Bott-Chern \'{e}quivariants singuliers de Bismut offre
un moyen unique de d\'{e}finir une th\'{e}orie des classes de
Bott-Chern \'{e}quivariantes singuli\`{e}res. Ceci
g\'{e}n\'{e}ralise la discussion de J. I. Burgos Gil et R.
Li\c{t}canu dans le cas \'{e}quivariant. En tant qu'un
sous-produit de cette \'{e}tude, nous allons prouver une formule
de concentration qui peut \^{e}tre utilis\'{e}e pour prouver un
th\'{e}or\`{e}me de concentration arithm\'{e}tique en
g\'{e}om\'{e}trie d'Arakelov.



\section{Introduction}
The Bott-Chern secondary characteristic classes were introduced by
R. Bott and S. S. Chern. They can be used to solve the problem
that the Chern-Weil theory is not additive for short exact
sequence of hermitian vector bundles. More precisely, assume that
we are given a short exact sequence
\begin{displaymath}
\overline{\varepsilon}:\quad
0\rightarrow\overline{E}'\rightarrow\overline{E}\rightarrow\overline{E}''\rightarrow
0
\end{displaymath}
of hermitian vector bundles on a compact complex manifold $X$.
Then the alternating sum of Chern character forms ${\rm
ch}(\overline{E}')-{\rm ch}(\overline{E})+{\rm
ch}(\overline{E}'')$ is not equal to $0$ unless this sequence is
orthogonally split. A Bott-Chern secondary characteristic class
associated to $\overline{\varepsilon}$ is an element
$\widetilde{{\rm ch}}(\overline{\varepsilon})\in \widetilde{A}(X)$
(cf. Section 2) satisfying

(i). (Differential equation) ${\rm dd}^c\widetilde{{\rm
ch}}(\overline{\varepsilon})={\rm ch}(\overline{E}')-{\rm
ch}(\overline{E})+{\rm ch}(\overline{E}'')$. Here the symbol ${\rm
dd}^c$ is the differential operator
$\frac{\overline{\partial}\partial}{2\pi i}$.

J.-M. Bismut, H. Gillet and C. Soul\'{e}'s construction of
Bott-Chern secondary classes (cf. \cite{BGS1}) forces it to
satisfy other two properties

(ii). (Functoriality) $f^*\widetilde{{\rm
ch}}(\overline{\varepsilon})=\widetilde{{\rm
ch}}(f^*\overline{\varepsilon})$ if $f: X'\rightarrow X$ is a
holomorphic map of complex manifolds.

(iii). (Normalization) $\widetilde{{\rm
ch}}(\overline{\varepsilon})=0$ if $\overline{\varepsilon}$ is
orthogonally split.

It has been shown that the three properties above are already
enough to characterize a theory of Bott-Chern secondary
characteristic classes. The same thing goes to the Chern-Weil
theory in the equivariant case. We shall recall these results in
Section 2 for the convenience of the reader.

In \cite{BGS3}, J.-M. Bismut, H. Gillet and C. Soul\'{e} defined
the Bott-Chern singular currents in order to solve a similar
differential equation as in (i) with respect to the resolution of
hermitian vector bundle associated to a closed immersion of
complex manifolds. Later in \cite{Bi}, J.-M. Bismut generalized
this topic to the equivariant case. Precisely speaking, let $G$ be
a compact Lie group and let $i: Y\rightarrow X$ be a
$G$-equivariant closed immersion of complex manifolds with
hermitian normal bundle $\overline{N}$. Suppose that
$\overline{\eta}$ is an equivariant hermitian vector bundle on $Y$
and that $\overline{\xi}.$ is a complex of equivariant hermitian
vector bundles providing a resolution of $i_*\overline{\eta}$ on
$X$ whose metrics satisfy Bismut assumption (A). Then, fixing an
element $g\in G$, J.-M. Bismut can construct a singular current
$T_g(\overline{\xi}.)\in D(X_g)$ which is a sum of $(p,p)$-type
currents satisfying

(i'). (Differential equation) ${\rm
dd}^cT_g(\overline{\xi}.)={i_g}_*({\rm ch}_g(\overline{\eta}){\rm
Td}_g^{-1}(\overline{N}))-\sum_k(-1)^k{\rm
ch}_g(\overline{\xi}_k)$.

As in the case of Bott-Chern secondary characteristic classes, it
can be shown that $T_g(\overline{\xi}.)$ also satisfies other two
properties

(ii'). (Functoriality)
$f_g^*T_g(\overline{\xi}.)=T_g(f^*\overline{\xi}.)$ if $f:
X'\rightarrow X$ is a $G$-equivariant holomorphic map of complex
manifolds which is transversal to $Y$.

(iii') (Normalization) $T_g(\overline{\xi.})=-\widetilde{{\rm
ch}}_g(\overline{\xi.})$ if $Y$ is the empty set.

Naturally, one hopes that such three properties are enough to
characterize a theory of equivariant singular Bott-Chern classes.
But unfortunately this is not true, $T_g(\overline{\xi}.)$ is not
the unique element which satisfies the properties (i'), (ii') and
(iii') even in the current class space
$\widetilde{\mathcal{U}}(X_g)$.

For the non-equivariant case i.e. when $G$ is the trivial group,
J. I. Burgos Gil and R. Li\c{t}canu have obtained a satisfactory
axiomatic characterization of singular Bott-Chern classes in their
article \cite{BL}. They realized this by adding a natural fourth
axiom to the properties (i'), (ii') and (iii') (removing the
subscript $g$) which is called the condition of homogeneity. In
this paper, we will do the equivariant version. Our strategy is
basically the same as that was used in J. I. Burgos Gil and R.
Li\c{t}canu's article. By deforming a resolution to a easily
understandable one, we show that a theory of equivariant singular
Bott-Chern classes is totally determined by its effects on Koszul
resolutions (cf. Theorem~\ref{601}). This approach can be viewed
as an analogue of J.-M. Bismut, H. Gillet and C. Soul\'{e}'s
axiomatic construction of Bott-Chern secondary characteristic
classes.

The original purpose of the author's study of the uniqueness
property of equivariant singular Bott-Chern classes is that he
wants to prove a purely analytic statement which is called the
concentration formula. Such a formula plays a crucial role in the
proof of an arithmetic concentration theorem in Arakelov geometry.
We shall formulate this result in the last section of this paper.

\textbf{Acknowledgements.} The author wishes to thank Damian
Roessler and Xiaonan Ma who suggested him to pay attention to
Jos\'{e} I. Burgos Gil and R\v{a}zvan Li\c{t}canu's work on the
uniqueness problem of singular Bott-Chern classes. The author is
also grateful to Jos\'{e} I. Burgos Gil for many fruitful
discussions between them, for his careful reading of a early
version of this paper, and for his suggestions which improve the
results of a crucial lemma.

\section{Equivariant secondary characteristic classes}
To every hermitian vector bundle on a compact complex manifold we
can associate a smooth differential form by using Chern-Weil
theory. Notice that Chern-Weil theory is not additive for short
exact sequence of hermitian vector bundles, the Bott-Chern
secondary characteristic classes cover this gap. In this section,
we shall recall how to generalize all these things above to the
equivariant setting, namely for a compact complex manifold $X$
which admits a holomorphic action of a compact Lie group $G$.

Let $g\in G$ be an automorphism of $X$, we denote by $X_g=\{x\in
X\mid g\cdot x=x\}$ the fixed point submanifold. $X_g$ is also a
compact complex manifold. Let $\overline{E}$ be an equivariant
hermitian vector bundle on $X$, this means that $E$ admits a
$G$-action which is compatible with the $G$-action on $X$ and that
the metric on $E$ is invariant under the action of $G$. If there
is no additional description, a morphism between equivariant
vector bundles will be a morphism of vector bundles which is
compatible with the equivariant structures. It is well known that
the restriction of an equivariant hermitian vector bundle
$\overline{E}$ to $X_g$ splits as a direct sum
\begin{displaymath}
\overline{E}\mid_{X_g}=\bigoplus_{\zeta\in
S^1}\overline{E}_{\zeta}
\end{displaymath}
where the equivariant structure $g^E$ of $E$ acts on
$\overline{E}_{\zeta}$ as $\zeta$. We often write $\overline{E}_g$
for $\overline{E}_1$ and denote its orthogonal complement by
$\overline{E}_{\bot}$. As usual, $A^{p,q}(X)$ stands for the space
of $(p,q)$-forms
$\Gamma^{\infty}(X,\Lambda^pT^{*(1,0)}X\wedge\Lambda^qT^{*(0,1)}X)$,
we define
\begin{displaymath}
\widetilde{A}(X)=\bigoplus_{p=0}^{{\rm dim}X}(A^{p,p}(X)/({\rm
Im}\partial+{\rm Im}\overline{\partial})).
\end{displaymath}
We denote by $\Omega^{\overline{E}_\zeta}$ the curvature matrix
associated to $\overline{E}_\zeta$. Let $(\phi_\zeta)_{\zeta\in
S^1}$ be a family of $\GL(\C)$-invariant formal power series such
that $\phi_\zeta\in \C[[\mathbf{gl}_{{\rm rk}E_\zeta}(\C)]]$ where
${\rm rk}E_\zeta$ stands for the rank of $E_\zeta$ which is a
locally constant function on $X_g$. Moreover, let $\phi\in
\C[[\bigoplus_{\zeta\in S^1}\C]]$ be any formal power series. We
have the following definition.

\begin{defn}\label{201}
The way to associate a smooth form to an equivariant hermitian
vector bundle $\overline{E}$ by setting
\begin{displaymath}
\phi_g(\overline{E}):=\phi((\phi_{\zeta}(-\frac{\Omega^{\overline{E}_\zeta}}{2\pi
i}))_{\zeta\in S^1})
\end{displaymath}
is called an equivariant
Chern-Weil theory associated to $(\phi_\zeta)_{\zeta\in S^1}$ and
$\phi$. The class of $\phi_g(\overline{E})$ in
$\widetilde{A}(X_g)$ is independent of the metric.
\end{defn}

The theory of equivariant secondary characteristic classes is
described in the following theorem.

\begin{thm}\label{202}
To every short exact sequence $\overline{\varepsilon}:
0\rightarrow \overline{E}'\rightarrow \overline{E}\rightarrow
\overline{E}''\rightarrow 0$ of equivariant hermitian vector
bundles on $X$, there is a unique way to attach a class
$\widetilde{\phi}_g(\overline{\varepsilon})\in \widetilde{A}(X_g)$
which satisfies the following three conditions:

(i). $\widetilde{\phi}_g(\overline{\varepsilon})$ satisfies the
differential equation
\begin{displaymath}
{\rm
dd}^c\widetilde{\phi}_g(\overline{\varepsilon})=\phi_g(\overline{E}'\oplus\overline{E}'')-\phi_g(\overline{E});
\end{displaymath}

(ii). for every equivariant holomorphic map $f: X'\rightarrow X$,
$\widetilde{\phi}_g(f^*\overline{\varepsilon})=f_g^*\widetilde{\phi}_g(\overline{\varepsilon})$;

(iii). $\widetilde{\phi}_g(\overline{\varepsilon})=0$ if
$\overline{\varepsilon}$ is equivariantly and orthogonally split.
\end{thm}
\begin{proof}
Firstly note that one can carry out the principle of \cite[Section
f.]{BGS1} to construct a new exact sequence of equivariant
hermitian vector bundles
\begin{displaymath}
\overline{\widetilde{\varepsilon}}:\quad 0\rightarrow
\overline{E'(1)}\rightarrow \overline{\widetilde{E}}\rightarrow
\overline{E}''\rightarrow 0
\end{displaymath} on $X\times
\mathbb{P}^1$ such that $i_0^*\overline{\widetilde{\varepsilon}}$
is isometric to $\overline{\varepsilon}$ and
$i_\infty^*\overline{\widetilde{\varepsilon}}$ is equivariantly
and orthogonally split. Here the projective line $\mathbb{P}^1$
carries the trivial $G$-action and the section $i_0$ (resp.
$i_\infty$) is defined by setting $i_0(x)=(x,0)$ (resp.
$i_\infty(x)=(x,\infty)$). Then one can show that an equivariant
secondary characteristic class
$\widetilde{\phi}_g(\overline{\varepsilon})$ which satisfies the
three conditions in the statement of this theorem must be of the
form
\begin{displaymath}
\widetilde{\phi}_g(\overline{\varepsilon})=-\int_{\mathbb{P}^1}\phi_g(\widetilde{E},h^{\widetilde{E}})\cdot\log\mid
z\mid^2.
\end{displaymath} So the uniqueness has been proved. For the existence, one may take this identity as the
definition of the equivariant secondary class
$\widetilde{\phi}_g$, of course one should verify that this
definition is independent of the choice of the metric
$h^{\widetilde{E}}$ and really satisfies the three conditions
above. The verification is totally the same as the non-equivariant
case, one just add the subscript $g$ to every corresponding
notation.

Another way to show the existence is to use the non-equivariant
secondary classes on $X_g$ directly. We first split
$\overline{\varepsilon}$ on $X_g$ into a family of short exact
sequences
\begin{displaymath}
\overline{\varepsilon}_\zeta:\quad 0\rightarrow
\overline{E}'_\zeta\rightarrow \overline{E}_\zeta\rightarrow
\overline{E}''_\zeta\rightarrow 0 \end{displaymath} for all
$\zeta\in S^1$. Using the non-equivariant secondary classes on
$X_g$ we define for $\zeta, \eta\in S^1$
\begin{displaymath}
(\widetilde{\phi_\zeta+\phi_\eta})(\overline{\varepsilon}_\zeta,\overline{\varepsilon}_\eta):=
\widetilde{\phi}_\zeta(\overline{\varepsilon}_\zeta)+\widetilde{\phi}_\eta(\overline{\varepsilon}_\eta)
\end{displaymath} and
\begin{displaymath}
(\widetilde{\phi_\zeta\cdot\phi_\eta})(\overline{\varepsilon}_\zeta,\overline{\varepsilon}_\eta):=
\widetilde{\phi}_\zeta(\overline{\varepsilon}_\zeta)\cdot\phi_\eta(\overline{E}_\eta)+\phi_\zeta(\overline{E}'_\zeta
+\overline{E}''_\zeta)\cdot\widetilde{\phi}_\eta(\overline{\varepsilon}_\eta)
\end{displaymath} and similarly for other finite sums and products.
With these notations we define
$\widetilde{\phi}_g(\overline{\varepsilon}):=\widetilde{\phi((\phi_\zeta)_{\zeta\in
S^1})}((\overline{\varepsilon}_\zeta)_{\zeta\in S^1})$. The
equivariant secondary class $\widetilde{\phi}_g$ defined like this
way satisfies the three conditions in the statement of this
theorem, this fact follows from the axiomatic characterization of
non-equivariant secondary classes.
\end{proof}

\begin{rem}\label{203}
(i). The first way to construct equivariant secondary
characteristic classes is also valid for long exact sequences of
hermitian vector bundles $\overline{\varepsilon}: 0\rightarrow
\overline{E}_m\rightarrow \overline{E}_{m-1}\rightarrow
\cdots\rightarrow \overline{E}_1\rightarrow
\overline{E}_0\rightarrow 0$. Here the sign is chosen so that
\begin{displaymath}
{\rm
dd}^c\widetilde{\phi}_g(\overline{\varepsilon})=\phi_g(\bigoplus_{j\
even}\overline{E}_j)-\phi_g(\bigoplus_{j\ odd}\overline{E}_j).
\end{displaymath} That means there exists an exact sequence
$\overline{\widetilde{\varepsilon}}$ on $X\times \mathbb{P}^1$
such that $i_0^*\overline{\widetilde{\varepsilon}}$ is isometric
to $\overline{\varepsilon}$ and
$i_\infty^*\overline{\widetilde{\varepsilon}}$ is equivariantly
and orthogonally split. This new exact sequence is called the
first transgression exact sequence of $\overline{\varepsilon}$ and
will be denoted by ${\rm tr}_1(\overline{\varepsilon})$.

(ii). The first part of this remark gives a uniqueness theorem for
secondary classes for long exact sequences. Then when $\phi_g$ is
additive one can have another way to construct the secondary
classes, that is to split a long exact sequence into a series of
short exact sequences and use the secondary classes in
Theorem~\ref{202} to formulate an alternating sum. This
alternating sum provides a secondary class for original long exact
sequence.
\end{rem}

We now give some examples of equivariant character forms and their
corresponding secondary characteristic classes.

\begin{exa}\label{204}
(i). The equivariant Chern character form ${\rm
ch}_g(\overline{E}):=\sum_{\zeta\in S^1}\zeta{\rm
ch}(\overline{E}_\zeta)$.

(ii). The equivariant Todd form ${\rm
Td}_g(\overline{E}):=\frac{c_{{\rm rk}E_g}(\overline{E}_g)}{{\rm
ch}_g(\sum_{j=0}^{{\rm rk}E}(-1)^j\wedge^j\overline{E}^\vee)}$. As
in \cite[Thm. 10.1.1]{Hi} one can show that
\begin{displaymath}
{\rm Td}_g(\overline{E})={\rm
Td}(\overline{E}_g)\prod_{\zeta\neq1}{\rm
det}(\frac{1}{1-\zeta^{-1}e^{\frac{\Omega^{\overline{E}_\zeta}}{2\pi
i}}}). \end{displaymath}

(iii). Let $\overline{\varepsilon}: 0\rightarrow
\overline{E}'\rightarrow \overline{E}\rightarrow
\overline{E}''\rightarrow 0$ be a short exact sequence of
hermitian vector bundles. The secondary Bott-Chern characteristic
class is given by $\widetilde{{\rm
ch}}_g(\overline{\varepsilon})=\sum_{\zeta\in S^1}\zeta
\widetilde{{\rm ch}}(\overline{\varepsilon}_\zeta)$.

(iv). If the equivariant structure $g^\varepsilon$ has the
eigenvalues $\zeta_1,\cdots,\zeta_m$, then the secondary Todd
class is given by
\begin{displaymath}
\widetilde{{\rm
Td}}_g(\overline{\varepsilon})=\sum_{i=1}^m\prod_{j=1}^{i-1}{\rm
Td}_g(\overline{E}_{\zeta_j})\cdot\widetilde{{\rm
Td}}(\overline{\varepsilon}_{\zeta_i})\cdot\prod_{j=i+1}^m{\rm
Td}_g(\overline{E}'_{\zeta_j}+\overline{E}''_{\zeta_j}).
\end{displaymath}
\end{exa}

\begin{rem}\label{205}
One can use Theorem~\ref{202} to give a proof of the statements
(iii) and (iv) in the examples above.
\end{rem}

\begin{lem}\label{206}
Let $\overline{\varepsilon}$ be an acyclic complex of equivariant
hermitian vector bundles. Then for any non-negative integer $k$ we
have
\begin{displaymath}
\widetilde{\phi}_g(\overline{\varepsilon}[-k])=(-1)^k\widetilde{\phi}_g(\overline{\varepsilon}).
\end{displaymath}
\end{lem}
\begin{proof}
This follows from the construction of $\overline{\varepsilon}[-k]$
which is obtained by shifting degree.
\end{proof}

A particular secondary class when we consider a fixed vector
bundle with two different hermitian metrics will be used
frequently in our paper, so we formulate it separately in the
following definition.

\begin{defn}\label{207}
Let $E$ be an equivariant vector bundle on $X$. Assume that $h_0$
and $h_1$ are two invariant hermitian metrics on $E$. We denote by
$\widetilde{\phi}_g(E,h_0,h_1)$ the equivariant secondary
characteristic class associated to the short exact sequence
\begin{displaymath}
0\rightarrow 0\rightarrow (E,h_1)\rightarrow (E,h_0)\rightarrow 0
\end{displaymath} so that we have the differential equation ${\rm
dd}^c\widetilde{\phi}_g(E,h_0,h_1)=\phi_g(E,h_0)-\phi_g(E,h_1)$.
\end{defn}

The following proposition describes the additivity of equivariant
secondary characteristic classes.

\begin{prop}\label{208}
Let
\begin{displaymath}
\xymatrix{ & 0 \ar[d] & 0 \ar[d] & 0 \ar[d] & \\
0 \ar[r] & \overline{E}_1' \ar[r] \ar[d] & \overline{E}_1 \ar[r]
\ar[d] & \overline{E}_1'' \ar[r] \ar[d] & 0 \\
0 \ar[r] & \overline{E}_2' \ar[r] \ar[d] & \overline{E}_2 \ar[r]
\ar[d] & \overline{E}_2'' \ar[r] \ar[d] & 0 \\
0 \ar[r] & \overline{E}_3' \ar[r] \ar[d] & \overline{E}_3 \ar[r]
\ar[d] & \overline{E}_3'' \ar[r] \ar[d] & 0 \\
& 0 & 0 & 0 &} \end{displaymath} be a double complex of
equivariant hermitian vector bundles on $X$ where all rows
$\overline{\varepsilon}_i$ and all columns $\overline{\delta}_j$
are exact. Then we have
\begin{displaymath}
\widetilde{\phi}_g(\overline{\varepsilon}_1\oplus\overline{\varepsilon}_3)-\widetilde{\phi}_g(\overline{\varepsilon}_2)
=\widetilde{\phi}_g(\overline{\delta}_1\oplus\overline{\delta}_3)-\widetilde{\phi}_g(\overline{\delta}_2).
\end{displaymath}
\end{prop}
\begin{proof}
We may have the corresponding diagram of hermitian vector bundles
on $X\times \mathbb{P}^1$ by the first construction in the proof
of Theorem~\ref{202}. Then
\begin{align*}
\widetilde{\phi}_g(\overline{\varepsilon}_2)-\widetilde{\phi}_g(\overline{\varepsilon}_1\oplus\overline{\varepsilon}_3)
&=-\int_{\mathbb{P}^1}[\phi_g(\widetilde{E_2},h^{\widetilde{E_2}})-\phi_g(\widetilde{E_1}\oplus\widetilde{E_3},h^{\widetilde{E_1}}\oplus
h^{\widetilde{E_3}})]\cdot\log\mid
z\mid^2\\
&=\int_{\mathbb{P}^1}{\rm
dd}^c\widetilde{\phi}_g(\overline{\widetilde{\delta_2}})\cdot\log\mid
z\mid^2=\int_{\mathbb{P}^1}\widetilde{\phi}_g(\overline{\widetilde{\delta_2}})\cdot{\rm
dd}^c\log\mid
z\mid^2\\
&=i_0^*\widetilde{\phi}_g(\overline{\widetilde{\delta_2}})-i_\infty^*\widetilde{\phi}_g(\overline{\widetilde{\delta_2}})
=\widetilde{\phi}_g(\overline{\delta}_2)-\widetilde{\phi}_g(\overline{\delta}_1\oplus\overline{\delta}_3).
\end{align*}
\end{proof}

\begin{rem}\label{209}
This proposition can be generalized without any difficulty to the
case of a bounded exact sequence of bounded exact sequences of
equivariant hermitian vector bundles. Let $\overline{A}_{*,*}$ be
such a double acyclic complex, we have
\begin{displaymath}
\widetilde{\phi}_g(\bigoplus_{k\
even}\overline{A}_{k,*})-\widetilde{\phi}_g(\bigoplus_{k\
odd}\overline{A}_{k,*})=\widetilde{\phi}_g(\bigoplus_{k\
even}\overline{A}_{*,k})-\widetilde{\phi}_g(\bigoplus_{k\
odd}\overline{A}_{*,k}).
\end{displaymath}
\end{rem}

\begin{cor}\label{210}
Let $\overline{A}_{*,*}$ be a bounded double complex of
equivariant hermitian vector bundles with exact rows and exact
columns, then we have
\begin{displaymath}
\widetilde{\phi}_g({\rm
Tot}\overline{A}_{*,*})=\widetilde{\phi}_g(\bigoplus_k\overline{A}_{k,*}[-k]).
\end{displaymath}
\end{cor}
\begin{proof}
For any non-negative integer $n$ we denote by ${\rm Tot}_n={\rm
Tot}((\overline{A}_{k,*})_{k\geq n})$ the total complex of the
exact complex formed by the rows with index bigger than $n-1$.
Then ${\rm Tot}_0={\rm Tot}(\overline{A}_{*,*})$. By an argument
of induction, for each $k\geq0$ we have an exact sequence of
complexes
\begin{displaymath}
0\rightarrow {\rm Tot}_{k+1}\rightarrow {\rm
Tot}_k\oplus\bigoplus_{l<k}\overline{A}_{l,*}[-l]\rightarrow
\bigoplus_{l\leq k}\overline{A}_{l,*}[-l]\rightarrow 0
\end{displaymath} which is orthogonally split in each degree.
Therefore by Proposition~\ref{208} and Remark~\ref{209} we get
\begin{displaymath}
\widetilde{\phi}_g({\rm
Tot}_k\oplus\bigoplus_{l<k}\overline{A}_{l,*}[-l])=\widetilde{\phi}_g({\rm
Tot}_{k+1}\oplus\bigoplus_{l\leq k}\overline{A}_{l,*}[-l]).
\end{displaymath} By induction, when $k$ is taken to be big enough
we prove the statement.
\end{proof}

\section{Cohomology of currents with fixed wave front sets}
This section is devoted to recall the results of \cite[Section
4]{BL} and to derive some standard consequences. To be more
precise, we recall that there is a classic theorem concerning the
complex of currents on a compact complex manifold which says that
its cohomology groups are isomorphic to the cohomology groups of
the complex of smooth forms. In this section, we shall prove a
similar theorem for the currents with any fixed wave front set.
This theorem implies a certain
$\partial\overline{\partial}$-lemma.

Let $X$ be a compact complex manifold of dimension $d$. Then the
space $A^n(X)$ of $C^\infty$ complex valued $n$-forms on $X$ is a
topological vector space with Schwartz topology (cf. \cite[Chapter
IX]{deRh}). We denote by $D_n(X)$ the continuous dual of $A^n(X)$
which is called the space of currents of dimension $n$ on $X$.
Note that $X$ is a complex manifold, we have the following
decomposition
\begin{displaymath}
A^n(X)=\bigoplus_{p+q=n}A^{p,q}(X)
\end{displaymath} and Dolbeault
operators $\partial: A^{p,q}(X)\rightarrow A^{p+1,q}(X),
\overline{\partial}: A^{p,q}(X)\rightarrow A^{p,q+1}(X)$ with
${\rm d}=\partial+\overline{\partial}$ from $A^n(X)$ to
$A^{n+1}(X)$ the usual differentials.

All things above induce corresponding notations for $D_n(X)$ as
follows
\begin{displaymath}
D_n(X)=\bigoplus_{p+q=n}D_{p,q}(X)
\end{displaymath} and Dolbeault
operators $\partial': D_{p+1,q}(X)\rightarrow D_{p,q}(X),
\overline{\partial}': D_{p,q+1}(X)\rightarrow D_{p,q}(X)$ with
${\rm d}'=\partial'+\overline{\partial}'$ from $D_{n+1}(X)$ to
$D_n(X)$. Here the differential ${\rm d}'$ should be understood as
for any $T\in D_{n+1}(X), \alpha\in A^n(X)$, ${\rm
d}'T(\alpha)=T({\rm d}\alpha)$. We now give two basic examples of
currents which will be used frequently.

\begin{exa}\label{301}
If $i: Y\hookrightarrow X$ is a $k-$dimensional analytic subspace
of $X$, we may define a $2k-$dimensional current $\delta_Y$ which
is introduced by Lelong \cite{Le} by
\begin{displaymath}
\delta_Y(\alpha)=\int_{Y^{ns}}i^*\alpha,\quad \alpha\in A^{2k}(X)
\end{displaymath} where $Y^{ns}$ is the subset of non-singular
points in $Y$. Note that $\delta_Y$ actually belongs to
$D_{k,k}(X)$ since if $\alpha^{p,q}\in A^{p,q}(X)$ with $p+q=2k$,
then $i^*\alpha=0$ unless $p=q=k$.
\end{exa}

\begin{exa}\label{302}
We may have the following products
\begin{displaymath}
D_n(X)\otimes A^m(X)\rightarrow D_{n-m}(X)
\end{displaymath} which
decompose into
\begin{displaymath}
D_{p,q}(X)\otimes A^{r,s}(X)\rightarrow D_{p-r,q-s}(X).
\end{displaymath} Actually for $T\in D_n(X), \alpha\in A^m(X)$, we
denote their product by $T\wedge \alpha$, and if $\beta\in
A^{n-m}(X)$, the product is defined by
\begin{displaymath}
(T\wedge\alpha)(\beta)=T(\alpha\wedge\beta).
\end{displaymath} In
particular, we have a map from $A^{p,q}(X)$ to $D_{d-p,d-q}(X)$
which maps $\alpha$ to $\delta_X\wedge\alpha$. We often write
$\delta_X\wedge\alpha$ as $[\alpha]$ for simplicity. From
\cite[Chapter X]{deRh} we know that the spaces $D_{p,q}$ have a
natural topology, for which the maps $A^{p,q}(X)\rightarrow
D_{d-p,d-q}(X)$ are continuous with dense images. So if we write
$D^{p,q}(X)=D_{d-p,d-q}(X)$, we may have the following embedding
\begin{displaymath}
A^{p,q}(X)\hookrightarrow D^{p,q}(X).
\end{displaymath} We would
like to indicate that, more generally, if $\alpha$ is a $L^1$-form
i.e. $\alpha$ has coefficients which are locally integrable then
$[\alpha]$ is a well-defined current.
\end{exa}

\begin{rem}\label{303}
The map $\alpha\mapsto[\alpha]$ doesn't send ${\rm d}$ to ${\rm
d}'$. In fact, for $\alpha\in A^n(X)$ and $\beta\in A^{d-n-1}(X)$,
by Stokes theorem we will have
\begin{align*}
[{\rm d}\alpha](\beta)&=\int_X{\rm d}\alpha\wedge\beta=\int_X{\rm
d}(\alpha\wedge\beta)-\int_X(-1)^n\alpha\wedge {\rm
d}\beta\\
&=(-1)^{n+1}\int_X\alpha\wedge {\rm d}\beta=(-1)^{n+1}({\rm
d}'[\alpha])(\beta). \end{align*} So if we write ${\rm
d}=(-1)^{n+1}{\rm d}'$ the differential from $D^n(X)$ to
$D^{n+1}(X)$, then the inclusion $A^n(X)\hookrightarrow D^n(X)$
commutes with ${\rm d}$. The same conclusions can be obtained for
$\partial$ and $\overline{\partial}$. And one should notice that
this commutativity induces a family of morphisms at the level of
cohomology with respect to $\partial,\overline{\partial}$ and
${\rm d}$. These morphisms are actually isomorphisms.
\end{rem}

The wave front set ${\rm WF}(\eta)$ of a current $\eta$ is a
closed conical subset of $T_\R^*X_0:=T_\R^*X\backslash\{0\}$, the
real cotangent bundle removed the complete zero section. This
conical subset measures the singularities of $\eta$, actually the
projection of ${\rm WF}(\eta)$ in $X$ is equal to the singular
locus of the support of $\eta$. It also allows us to define
certain products and pull-backs of currents. Let $S$ be a conical
subset of $T_\R^*X_0$ and let $D^*(X,S)$ stand for the spaces
consisting of all currents whose wave front sets are contained in
$S$. Now suppose that $P$ is a differential operator with smooth
coefficients, then we have ${\rm WF}(P\circ\eta)\subseteq {\rm
WF}(\eta)$ by \cite[(8.1.11)]{Hoe}. This means $D^*(X,S)$ form a
$\partial$-, $\overline{\partial}$- or ${\rm d}$-complex.

Let $f: Y\rightarrow X$ be a morphism of compact complex
manifolds. The set of normal directions of $f$ is
\begin{displaymath}
N_f=\{(f(y),v)\in T_\R^*X \mid {\rm d}f^tv=0\}.
\end{displaymath}
This set measures the singularities of the morphism $f$. Actually,
if $f$ is smooth then $N_f=0$ and if $f$ is a closed immersion
then $N_f$ is the conormal bundle $N_{X/Y,\R}^\vee$. Let $S\subset
T_\R^*X_0$ be a closed conical subset, the morphism $f$ is
transversal to $S$ if $N_f\cap S=\emptyset$.

\begin{thm}\label{304}
Let $f: Y\rightarrow X$ be a morphism of compact complex manifolds
which is transversal to $S$. Then there is a unique way to extend
the pull-back $f^*: A^*(X)\rightarrow A^*(Y)$ to a continuous
morphism of complexes
\begin{displaymath}
f^*:\quad D^*(X,S)\rightarrow D^*(Y,f^*S).
\end{displaymath}
\end{thm}
\begin{proof}
This follows from \cite[Theorem 8.2.4]{Hoe}. Here the topology on
current space $D^*$ is given by \cite[Definition 8.2.2]{Hoe} which
is finer than the usual dual topology.
\end{proof}

\begin{thm}\label{306}($\overline{\partial}$-Poincar\'{e} lemma)
For any integer $p\geq 0$, denote by $\mathcal{D}_{X,S}^{p,*}$ the
sheaf of currents of type $(p,*)$ whose wave front sets are
contained in $S$. Then for each $q>0$, any
$\overline{\partial}$-closed section of $\mathcal{D}_{X,S}^{p,q}$
is locally $\overline{\partial}$-exact.
\end{thm}
\begin{proof}
Let $\mathcal{A}^{p,*}_X$ be the sheaf of smooth forms of type
$(p,*)$ on $X$, we claim that the natural inclusions
\begin{displaymath}
\iota: \mathcal{A}^{p,q}_X\rightarrow \mathcal{D}_{X,S}^{p,q}
\end{displaymath} induce a quasi-isomorphism between complexes.
Actually, this claim is the content of \cite[Theorem 4.5]{BL}.
With this observation we may reduce our problem to the case of
smooth forms which is classical, one can find in \cite[Page
25]{GH} a proof of this statement.
\end{proof}

\begin{cor}\label{307}
Let notations and assumptions be as in Theorem~\ref{306} and its
proof, then the natural morphisms
$H^{p,*}(\mathcal{A}^{p,*}_X(X),\overline{\partial})\longrightarrow
H^{p,*}(\mathcal{D}^{p,*}_{X,S}(X),\overline{\partial})$ are
isomorphisms.
\end{cor}
\begin{proof}
We denote by $\Omega^p$ the sheaf of holomorphic $p$-forms on $X$.
The Dolbeault theorem says that $H^*(X,\Omega^p)$ are isomorphic
to $H^{p,*}(\mathcal{A}^{p,*}_X(X),\overline{\partial})$. The
proof of the Dolbeault theorem relies on two deep results, one is
that the following complex of sheaves
\begin{displaymath}
\xymatrix{0 \ar[r] & \Omega^p \ar[r] & \mathcal{A}^{p,0}_X
\ar[r]^-{\overline{\partial}} & \mathcal{A}^{p,1}_X
\ar[r]^-{\overline{\partial}} & \cdots
\ar[r]^-{\overline{\partial}} & \mathcal{A}^{p,d}_X \ar[r] & 0}
\end{displaymath} is exact, the other one is that the sheaves
$\mathcal{A}^{p,*}_X$ admit partitions of unity so that
$H^k(X,\mathcal{A}^{p,*}_X)=0$ for $k>0$. Note that the sheaves
$\mathcal{D}^{p,*}_{X,S}$ may be multiplied by $C^\infty$
functions, hence they also admit partitions of unity. Therefore
one can carry out the principle of the sheaf-theoretic proof of
Dolbeault theorem to prove that $H^*(X,\Omega^p)\cong
H^{p,*}(\mathcal{D}^{p,*}_{X,S}(X),\overline{\partial})$ if the
following complex of sheaves
\begin{displaymath}
\xymatrix{0 \ar[r] & \Omega^p \ar[r] & \mathcal{D}^{p,0}_{X,S}
\ar[r]^-{\overline{\partial}} & \mathcal{D}^{p,1}_{X,S}
\ar[r]^-{\overline{\partial}} & \cdots
\ar[r]^-{\overline{\partial}} & \mathcal{D}^{p,d}_{X,S} \ar[r] &
0}
\end{displaymath} is exact. Such a Dolbeault theorem for
currents implies our statements in this corollary. Indeed, the
complex of sheaves above is really exact, the exactness at
$0$-degree is just the regularity theorem for the
$\overline{\partial}$-operator (cf. \cite[Page 380]{GH}) and the
exactness at higher degrees is implied by the
$\overline{\partial}$-Poincar\'{e} lemma, Theorem~\ref{306}.
\end{proof}

\begin{rem}\label{308}
One can prove the similar results for $\partial-$cohomology and de
Rham cohomology, namely the natural morphisms
$H^{*,p}(\mathcal{A}^{*,p}_X(X),\partial)\longrightarrow
H^{*,p}(\mathcal{D}^{*,p}_{X,S}(X),\partial)$ and $H^*_{\rm
DR}(X)\longrightarrow H^{*}(\mathcal{D}^{*}_{X,S}(X),{\rm d})$ are
all isomorphisms.
\end{rem}

\begin{cor}\label{309}
Let $\mathcal{D}^{p,*}_X$ be the sheaf of currents of type $(p,*)$
on $X$, then the natural morphisms
$H^{p,*}(\mathcal{D}^{p,*}_{X,S}(X),\overline{\partial})\longrightarrow
H^{p,*}(\mathcal{D}^{p,*}_X(X),\overline{\partial})$ are
isomorphisms.
\end{cor}
\begin{proof}
This follows from Corollary~\ref{307}.
\end{proof}

This Corollary implies the following
$\partial\overline{\partial}$-lemma.

\begin{thm}\label{310}
Let $X$ be a compact complex manifold and let $S$ be a closed
conical subset of $T_\R^*X_0$ . Then:

(i). If $\gamma$ is a current on $X$ such that
$\partial\overline{\partial}\gamma\in D^*(X,S)$, then there exist
currents $\alpha$ and $\beta$ such that
$\gamma=\omega+\partial\alpha+\overline{\partial}\beta$ with
$\omega\in D^*(X,S)$.

(ii). If $\omega$ is an element in $D^*(X,S)$ such that
$\omega=\partial u+\overline{\partial}v$ for currents $u$ and $v$,
then there exist currents $\alpha,\beta\in D^*(X,S)$ such that
$\omega=\partial\alpha+\overline{\partial}\beta$.
\end{thm}
\begin{proof}
(i). The hypothesis $\partial\overline{\partial}\gamma=\eta$ with
$\eta\in D^*(X,S)$ implies that
$\eta=\partial(\overline{\partial}\gamma)$ and hence
$\eta=\partial\alpha$ for some $\alpha\in D^*(X,S)$. So
$\partial(\overline{\partial}\gamma-\alpha)=0$ and
$\overline{\partial}\gamma-\alpha=\beta+\partial\gamma_1$ with
$\beta\in D^*(X,S)$. So we know that
$\partial\overline{\partial}\gamma_1=\eta_1=\overline{\partial}(\alpha+\beta)$
is contained in $D^*(X,S)$. By repeating this argument we get a
sequence of currents $\gamma_n$ such that
$\overline{\partial}\gamma_n=u_n+\partial\gamma_{n+1}$ with
$u_n\in D^*(X,S)$.

Note that if we assume that $\gamma_n\in D^{p,q}(X)$, then
$\gamma_{n+1}$ should be in $D^{p-1,q+1}(X)$ by construction. So
when $n$ is big enough we have $\gamma_{n+1}=0$. Therefore
$\overline{\partial}\gamma_n=u_n$ is contained in $D^*(X,S)$,
hence $\gamma_n=\omega_n+\overline{\partial}\beta_n$ with
$\omega_n\in D^*(X,S)$. So
$\overline{\partial}(\gamma_{n-1}+\partial\beta_n)=u_{n-1}+\partial\omega_n$
is contained in $D^*(X,S)$, and therefore
$\gamma_{n-1}=\omega_{n-1}+\partial\alpha_{n-1}+\overline{\partial}\beta_{n-1}$
with $\omega_{n-1}\in D^*(X,S)$. By repeating this argument we are
done.

(ii). If $\omega=\partial u+\overline{\partial}v$, then
$\partial\omega=\partial\overline{\partial}v$ which implies that
$v=\alpha+\partial x+\overline{\partial}y$ with $\alpha\in
D^*(X,S)$ by (i). So we have
$\overline{\partial}v=\overline{\partial}\alpha+\overline{\partial}\partial
x$. Similarly $\partial
u=\partial\beta+\partial\overline{\partial}z$ with $\beta\in
D^*(X,S)$. Therefore
$\omega=\partial\alpha+\overline{\partial}\beta+\partial\overline{\partial}(z-x)$.
Again by (i), $z-x=\gamma+\partial s+\overline{\partial}t$ with
$\gamma\in D^*(X,S)$. So
$\partial\overline{\partial}(z-x)=\partial\overline{\partial}\gamma$
which implies that
$\omega=\partial(\alpha+\overline{\partial}\gamma)+\overline{\partial}\beta$.
\end{proof}

\section{Deformation to the normal cone}
By a projective manifold we shall understand a compact complex
manifold which is projective algebraic, that means a projective
manifold is the complex analytic space $X(\C)$ associated to a
smooth projective variety $X$ over $\C$. Denote by $\mu_n$ the
diagonalisable group variety over $\C$ associated to $\Z/{n\Z}$,
we say $X$ is equivariant if it admits a $\mu_n$-projective action
(cf. \cite[Section 2]{KR1}). Write $X_{\mu_n}$ for the fixed point
subscheme, by GAGA principle, $X_{\mu_n}(\C)$ is equal to
$X(\C)_g$ where $g$ is the automorphism on $X(\C)$ corresponding
to a fixed primitive $n$-th root of unity. From now on, if no
confusion arises, we shall not distinguish between $X$ and $X(\C)$
as well as $X_{\mu_n}$ and $X_g$.

In this section, we shall describe the algebro-geometric
preliminaries for the discussion of the uniqueness of equivariant
singular Bott-Chern classes. Our main tool is an elegant method so
called the deformation to the normal cone which allows us to
deform a resolution of hermitian vector bundle associated to a
closed immersion of projective manifolds to a simpler one. This
will help us to formulate the analytic data (e.g. the secondary
characteristic class) of the original resolution by using the
corresponding analytic data of the new one. This process is just
like the first construction we mentioned in the proof of
Theorem~\ref{202}.

The first part of this section is devoted to recall the
deformation to the normal cone technique which can be found in
several standard literatures, for example in \cite[Section
4]{BGS2}. The second part is devoted to the equivariant analogue.

Let $i: Y\hookrightarrow X$ be a closed immersion of projective
manifolds. We will denote by $N_{X/Y}$ the normal bundle of this
immersion. For a vector bundle $E$ on $X$ or $Y$, the notation
$\mathbb{P}(E)$ will stand for the projective space bundle ${\rm
Proj}({\rm Sym}(E^\vee))$.

\begin{defn}\label{401}
The deformation to the normal cone $W(i)$ of the immersion $i$ is
the blowing up of $X\times \mathbb{P}^1$ along $Y\times
\{\infty\}$. We shall just write $W$ for $W(i)$ if there is no
confusion about the immersion.
\end{defn}

We denote by $p_X$ (resp. $p_Y$) the projection $X\times
\mathbb{P}^1\rightarrow X$ (resp. $Y\times \mathbb{P}^1\rightarrow
Y$) and by $\pi$ the blow-down map $W\rightarrow
X\times\mathbb{P}^1$. We also denote by $q_X$ (resp. $q_Y$) the
projection $X\times \mathbb{P}^1\rightarrow \mathbb{P}^1$ (resp.
$Y\times \mathbb{P}^1\rightarrow \mathbb{P}^1$) and by $q_W$ the
composition $q_X\circ\pi$. It is well known that the map $q_W$ is
flat and for $t\in \mathbb{P}^1$, we have
\begin{displaymath}
q_W^{-1}(t)\cong\left\{
\begin{array}{ll}
    X\times \{t\}, & \hbox{if $t\neq\infty$,} \\
    P\cup \widetilde{X}, & \hbox{if $t=\infty$,} \\
\end{array}%
\right.
\end{displaymath} where $\widetilde{X}$ is isomorphic to
the blowing up of $X$ along $Y$ and $P$ is isomorphic to the
projective completion of $N_{X/Y}$ i.e. the projective space
bundle $\mathbb{P}(N_{X/Y}\oplus\mathcal{O}_Y)$. Denote the
canonical projection from $\mathbb{P}(N_{X/Y}\oplus\mathcal{O}_Y)$
to $Y$ by $\pi_P$, then the morphism $\mathcal{O}_Y\rightarrow
N_{X/Y}\oplus\mathcal{O}_Y$ induces a canonical section $i_\infty:
Y\hookrightarrow \mathbb{P}(N_{X/Y}\oplus\mathcal{O}_Y)$ which is
called the zero section embedding. Moreover, let $j: Y\times
\mathbb{P}^1\rightarrow W$ be the canonical closed immersion
induced by $i\times {\rm Id}$, then the component $\widetilde{X}$
doesn't meet $j(Y\times \mathbb{P}^1)$ and the intersection of
$j(Y\times \mathbb{P}^1)$ and $P$ is exactly the image of $Y$
under the section $i_\infty$.

On $P=\mathbb{P}(N_{X/Y}\oplus\mathcal{O}_Y)$, there exists a
tautological exact sequence
\begin{displaymath}
0\rightarrow \mathcal{O}(-1)\rightarrow
\pi_P^*(N_{X/Y}\oplus\mathcal{O}_Y)\rightarrow Q\rightarrow 0
\end{displaymath} where $Q$ is the tautological quotient bundle.
This exact sequence and the inclusion
$\mathcal{O}_P\rightarrow\pi_P^*(N_{X/Y}\oplus\mathcal{O}_Y)$
induce a section $\sigma: \mathcal{O}_P\rightarrow Q$ which
vanishes along the zero section $i_\infty(Y)$. By duality we get a
morphism $Q^\vee\rightarrow \mathcal{O}_P$, and this morphism
induces the following exact sequence
\begin{displaymath}
0\rightarrow
\wedge^nQ^\vee\rightarrow\cdots\rightarrow\wedge^2Q^\vee\rightarrow
Q^\vee\rightarrow\mathcal{O}_P\rightarrow
{i_\infty}_*\mathcal{O}_Y\rightarrow 0
\end{displaymath} where $n$
is the rank of $Q$. Note that $i_\infty$ is a section of $\pi_P$
i.e. $\pi_P\circ i_\infty={\rm Id}$, the projection formula
implies the following definition.

\begin{defn}\label{402}
For any vector bundle $F$ on $Y$, the following complex of vector
bundles
\begin{displaymath}
0\rightarrow
\wedge^nQ^\vee\otimes\pi_P^*F\rightarrow\cdots\rightarrow\wedge^2Q^\vee\otimes\pi_P^*F\rightarrow
Q^\vee\otimes\pi_P^*F\rightarrow\pi_P^*F\rightarrow 0
\end{displaymath} provides a resolution of ${i_\infty}_*F$ on $P$.
This complex is called the Koszul resolution of ${i_\infty}_*F$
and will be denoted by $K(F,N_{X/Y})$. If the normal bundle
$N_{X/Y}$ admits some hermitian metric, then the tautological
exact sequence induces a hermitian metric on $Q$. If, moreover,
the bundle $F$ also admits a hermitian metric, then the Koszul
resolution is a complex of hermitian vector bundles and will be
denoted by $K(\overline{F},\overline{N}_{X/Y})$.
\end{defn}

We now summarize the most important result about the application
of the deformation to the normal cone.

\begin{thm}\label{403}
Let $i: Y\hookrightarrow X$ be a closed immersion of projective
manifolds, and let $W=W(i)$ be the deformation to the normal cone
of $i$. Assume that $\overline{\eta}$ is a hermitian vector bundle
on $Y$ and $\overline{\xi}.$ is a complex of hermitian vector
bundles which provides a resolution of $i_*\overline{\eta}$ on
$X$. Then there exists a complex of hermitian vector bundles ${\rm
tr}_1(\overline{\xi}.)$ on $W$ such that

(i). ${\rm tr}_1(\overline{\xi}.)$ provides a resolution of
$j_*p_Y^*(\overline{\eta})$ on $W$;

(ii). ${\rm tr}_1(\overline{\xi}.)\mid_{X\times\{0\}}$ is
isometric to the original complex $\overline{\xi}.$;

(iii). the restriction of ${\rm tr}_1(\overline{\xi}.)$ to
$\widetilde{X}$ is orthogonally split;

(iv). the restriction of ${\rm tr}_1(\overline{\xi}.)$ to $P$ fits
an exact sequence of resolutions on $P$
\begin{displaymath}
\xymatrix{0 \ar[r] &
\overline{A}. \ar[r] \ar[d] & {\rm tr}_1(\overline{\xi}.)\mid_P
\ar[r] \ar[d] &
K(\overline{\eta},\overline{N}_{X/Y}) \ar[r] \ar[d] & 0 \\
& 0 \ar[r] & {i_\infty}_*\overline{\eta} \ar[r]^-{=} &
{i_\infty}_*(\overline{\eta}) & }
\end{displaymath} where
$\overline{A}.$ is orthogonally split and
$K(\overline{\eta},\overline{N}_{X/Y})$ is the hermitian Koszul
resolution;

(v). when $Y=\emptyset$, ${\rm tr}_1(\overline{\xi}.)$ is the
first transgression exact sequence introduced in Remark~\ref{203};

(vi). Let $f: X'\rightarrow X$ be a morphism of projective
manifolds which is smooth or transversal to $Y$. Formulate the
following
Cartesian square
\begin{displaymath}
\xymatrix{Y' \ar[r]^-{i'} \ar[d] & X' \ar[d] \\
Y \ar[r]^-{i} & X}
\end{displaymath} and denote by $f_W$ the
induced morphism from $W'=W(i')$ to $W$, then we have
\begin{displaymath}
f_W^*({\rm tr}_1(\overline{\xi}.))={\rm tr}_1(f^*\overline{\xi}.).
\end{displaymath}
\end{thm}
\begin{proof}
If $E$ is a vector bundle on $X$, we shall denote by $E(i)$ the
vector bundle on $X\times \mathbb{P}^1$ given by
$E(i)=p_X^*E\otimes q_X^*\mathcal{O}_{\mathbb{P}^1}(i)$. Now let
$\widetilde{C}.$ be the complex of vector bundles on $X\times
\mathbb{P}^1$ given by
$\widetilde{C}_i=\xi_i(i)\oplus\xi_{i-1}(i-1)$ with differential
${\rm d}(a,b)=(b,0)$. Let $y$ be a section of
$\mathcal{O}_{\mathbb{P}^1}(1)$ vanishing only at infinity, then
on $X\times(\mathbb{P}^1\backslash \{\infty\})$ we may construct a
family of inclusions of vector bundles $\gamma_i:
\xi_i\hookrightarrow \widetilde{C}_i$ given by $s\mapsto (s\otimes
y^i,(-1)^i{\rm d}s\otimes y^{i-1})$.

On the other hand, define a complex of vector bundles
$\widetilde{D}.$ on $X\times \mathbb{P}^1$ by
$\widetilde{D}_i=\xi_{i-1}(i)\oplus \xi_{i-2}(i-1)$. The morphism
of complexes $\varphi: \widetilde{C}.\rightarrow\widetilde{D}.$
given by $\varphi(s,t)=({\rm d}s+(-1)^it\otimes y,{\rm d}t)$
induces a morphism of complexes on $W$
\begin{displaymath}
\phi: \pi^*\widetilde{C}.\longrightarrow \pi^*\widetilde{D}.
\end{displaymath} where $\pi$ is the blow-down map. Then ${\rm
tr}_1(\xi.)$ is defined as the kernel of $\phi$.

Over $\pi^{-1}(X\times(\mathbb{P}^1\backslash\{\infty\}))$ we
shall endow ${\rm tr}_1(\xi.)_i$ with the metric induced by the
identification with $\xi_i$. And over
$\pi^{-1}(X\times(\mathbb{P}^1\backslash\{0\}))$ we shall endow
${\rm tr}_1(\xi.)_i$ with the metric induced by $\widetilde{C}_i$.
Finally we glue together these two metrics by a partition of unity
so that we get a hermitian metric on ${\rm tr}_1(\xi.)$. We refer
to \cite[Section 4]{BGS2} for the proof of the statement that the
complex of hermitian vector bundles ${\rm tr}_1(\overline{\xi}.)$
constructed in the way above really satisfies those conditions in
our theorem.
\end{proof}

\begin{rem}\label{404}
(i). Assume that $X$ is a $\mu_n$-equivariant projective manifold
and $E$ is an equivariant locally free sheaf on $X$. Then
according to \cite[(1.4) and (1.5)]{Ko}, $\mathbb{P}(E)$ admits a
canonical $\mu_n$-equivariant structure such that the projection
map $\mathbb{P}(E)\rightarrow X$ is equivariant and the canonical
bundle $\mathcal{O}(1)$ admits an equivariant structure. Moreover,
let $Y\to X$ be an equivariant closed immersion of projective
manifolds, according to \cite[(1.6)]{Ko} the action of $\mu_n$ on
$X$ can be extended to the blowing up ${\rm Bl}_YX$ such that the
blow-down map is equivariant and the canonical bundle
$\mathcal{O}(1)$ admits an equivariant structure. So the
constructions of blowing up and the deformation to the normal cone
are both compatible with the equivariant setting.

(ii). Furthermore, by endowing $\mathbb{P}^1$ with the trivial
action, we would like to reformulate all results in this section
especially Theorem~\ref{403} in the equivariant setting. We first
claim that the constructions of all vector bundles and bundle
morphisms in Theorem~\ref{403} also fit the equivariant setting.
This follows from the fact that they are all constructed
canonically. For more details, see \cite[Exp. VII, Lemme 2.4,
Proposition 2.5 and Lemme 3.2]{GBI} as well as \cite[Lemma 4.1,
Remark (ii) p. 314 and (4.7) p. 315]{BGS2}. But unfortunately, the
local uniqueness of resolutions (cf. \cite[Theorem 8]{Ei}) may not
be valid for the equivariant case so that the local method used in
the proof of the statement that the restriction of ${\rm
tr}_1(\overline{\xi}.)$ to $\widetilde{X}$ is orthogonally split
is not compatible with the equivariant setting. We have to
formulate relative results and proofs in a different way.
\end{rem}

\begin{lem}\label{tt1}
Let $X$ be a $\mu_n$-equivariant projective manifold, then the
category of coherent $\mu_n$-modules on $X$ is an abelian
category. A complex of $\mu_n$-equivariant coherent sheaves on $X$
is exact if and only if the underlying complex of
$\mathcal{O}_X$-modules is exact.
\end{lem}
\begin{proof}
This follows from \cite[Lemma 1.3]{Ko}.
\end{proof}

\begin{lem}\label{tt2}
Let $X$ be a $\mu_n$-equivariant projective manifold. In other
words, $X$ is a projective manifold which admits an automorphism
$g$ of order $n$. Assume that
\begin{displaymath}
\overline{\varepsilon}: 0 \rightarrow \overline{L}\rightarrow
\overline{E}\rightarrow \overline{F}\rightarrow 0
\end{displaymath} is a short exact sequence of equivariant hermitian vector
bundles on $X$. If the underlying sequence of hermitian vector
bundles is orthogonally split, then $\overline{\varepsilon}$ is
equivariantly and orthogonally split on $X$.
\end{lem}

\begin{proof}
Denote by $f$ the bundle morphism $E\rightarrow F$, by assumption
$f$ is equivariant. Since the underlying sequence of hermitian
vector bundles is orthogonally split, there exists a bundle
morphism $h$ from $F$ to $E$ such that $f\circ h={\rm Id}_{F}$ and
$\overline{F}$ is isometric to its image under this morphism $h$.
We recall that the $g$-structure on $\overline{E}$ (resp.
$\overline{F}$) is an isometry $\sigma_E: g^*\overline{E}\to
\overline{E}$ (resp. $\sigma_F: g^*\overline{F}\to \overline{F}$)
which satisfies certain associativity properties. We define a
$g$-action on the morphisms of equivariant bundles as follows. Let
$u: M\to N$ be a morphism of equivariant bundles, then
\begin{displaymath}
g\bullet u:=\sigma_N\circ g^*u\circ\sigma_M^{-1}
\end{displaymath}
which is still a morphism from $M$ to $N$. By definition, $u$ is
equivariant if and only if $g\bullet u=u$. One can easily check
that $g\bullet(g\bullet u)=g^2\bullet u$. Now since the morphisms
$f$ and ${\rm Id}_{F}$ are both equivariant, we compute
\begin{align*}
{\rm Id}_F&=g\bullet {\rm Id}_F=g\bullet(f\circ h)\\
&=\sigma_F\circ g^*(f\circ h)\circ \sigma_F^{-1}=\sigma_F\circ g^*f\circ g^*h\circ \sigma_F^{-1} \\
&=\sigma_F\circ g^*f\circ\sigma_E^{-1}\circ\sigma_E\circ g^*h\circ
\sigma_F^{-1}=f\circ(g\bullet h).
\end{align*}
Replacing $g$ by $g^k$ from $k=2$ to $k=n$, we get a meaningful
average of $h$ and it satisfies the following identity
\begin{displaymath}
f\circ(\frac{\sum_{k=0}^{n-1}g^k\bullet h}{n})={\rm Id}_F.
\end{displaymath}
Therefore $\frac{1}{n}\sum_{k=0}^{n-1}g^k\bullet h$ is an
equivariant section of $f$ which still makes $\overline{F}$
isometric to its image, so we are done.
\end{proof}

\begin{rem}\label{tt3}
In general, if the action on $X$ is not of finite order or the
base field of $X$ has characteristic dividing $n$ then the proof
given for Lemma~\ref{tt2} fails. Nevertheless, we can show that
$\overline{\varepsilon}$ is always equivariantly and orthogonally
split on $X_g$.

Actually, the problem that $\overline{\varepsilon}$ may not be
equivariantly and orthogonally split on the whole manifold $X$
arises because $h$ may not be equivariant. Note that on the fixed
point submanifold $X_g$, the morphism $h\mid_{(F\mid_{X_g})}$ is
equivariant if and only if it maps $F_\zeta$ into $E_\zeta$ for
any $\zeta\in S^1$. But this is rather clear because $f$ is
equivariant and the restriction of $f\circ h$ on $F\mid_{X_g}$ is
exactly the identity map on $F\mid_{X_g}$. So we are done.
\end{rem}

Together with Lemma~\ref{tt1}, Lemma~\ref{tt2} and
Remark~\ref{404}, we have the following theorem which is an
analogue of Theorem~\ref{403} in the equivariant setting.

\begin{thm}\label{main}
Let $i: Y\hookrightarrow X$ be an equivariant closed immersion of
equivariant projective manifolds, and let $W=W(i)$ be the
deformation to the normal cone of $i$. Assume that
$\overline{\eta}$ is an equivariant hermitian vector bundle on $Y$
and $\overline{\xi}.$ is a complex of equivariant hermitian vector
bundles which provides a resolution of $i_*\overline{\eta}$ on
$X$. Then there exists a complex of equivariant hermitian vector
bundles ${\rm tr}_1(\overline{\xi}.)$ on $W$ such that

(i). ${\rm tr}_1(\overline{\xi}.)$ provides an equivariant
resolution of $j_*p_Y^*(\overline{\eta})$ on $W$;

(ii). ${\rm tr}_1(\overline{\xi}.)\mid_{X\times\{0\}}$ is
isometric to the original complex $\overline{\xi}.$;

(iii). the restriction of ${\rm tr}_1(\overline{\xi}.)$ to
$\widetilde{X}$ is equivariantly and orthogonally split;

(iv). the restriction of ${\rm tr}_1(\overline{\xi}.)$ to $P$ fits
an equivariant exact sequence of equivariant resolutions on $P$
\begin{displaymath}
\xymatrix{0 \ar[r] & \overline{A}. \ar[r] \ar[d] & {\rm
tr}_1(\overline{\xi}.)\mid_P \ar[r] \ar[d] &
K(\overline{\eta},\overline{N}_{X/Y}) \ar[r] \ar[d] & 0 \\
& 0 \ar[r] & {i_\infty}_*\overline{\eta} \ar[r]^-{=} &
{i_\infty}_*(\overline{\eta}) & }
\end{displaymath} where
$\overline{A}.$ is an equivariantly and orthogonally split
complex, $K(\overline{\eta},\overline{N}_{X/Y})$ is the hermitian
Koszul resolution;

(v). when $Y=\emptyset$, ${\rm tr}_1(\overline{\xi}.)$ is the
first transgression exact sequence introduced in Remark~\ref{203};

(vi). Let $f: X'\rightarrow X$ be an equivariant morphism of
equivariant projective manifolds which is smooth or transversal to
$Y$. Formulate the following
Cartesian square
\begin{displaymath}
\xymatrix{Y' \ar[r]^-{i'} \ar[d] & X' \ar[d] \\
Y \ar[r]^-{i} & X}
\end{displaymath} and denote by $f_W$ the
induced morphism from $W'=W(i')$ to $W$, then we have
\begin{displaymath}
f_W^*({\rm tr}_1(\overline{\xi}.))={\rm tr}_1(f^*\overline{\xi}.).
\end{displaymath}
\end{thm}

To end this section, we recall some basic facts concerning the
relation between equivariant setting and non-equivariant setting.
Their proofs can be found in \cite[Section 2 and 6.2]{KR1}.

\begin{prop}\label{405}
Let $i: Y\hookrightarrow X$ be an equivariant closed immersion of
projective manifolds, and let $i_g: Y_g\hookrightarrow X_g$ be the
induced closed immersion between fixed point submanifolds. Then we
have

(i). the natural morphism $N_{X_g/{Y_g}}\rightarrow (N_{X/Y})_g$
is an isomorphism;

(ii). the natural morphism from the deformation to the normal cone
$W(i_g)$ to the fixed point submanifold $W(i)_g$ is a closed
immersion, this closed immersion induces the closed immersions
$\mathbb{P}(N_{X_g/{Y_g}}\oplus
\mathcal{O}_{Y_g})\rightarrow\mathbb{P}(N_{X/Y}\oplus\mathcal{O}_Y)_g$
and $\widetilde{X_g}\rightarrow \widetilde{X}_g$;

(iii). the fixed point submanifold of
$\mathbb{P}(N_{X/Y}\oplus\mathcal{O}_Y)$ is
$\mathbb{P}(N_{X_g/{Y_g}}\oplus \mathcal{O}_{Y_g})\coprod_{\zeta
\neq1} \mathbb{P}((N_{X/Y})_\zeta)$;

(iv). the closed immersion $i_{\infty,g}$ factors through
$\mathbb{P}(N_{X_g/{Y_g}}\oplus \mathcal{O}_{Y_g})$ and the other
components $\mathbb{P}((N_{X/Y})_\zeta)$ don't meet $Y$. Hence the
complex $K(\mathcal{O}_Y,N_{X/Y})_g$, obtained by taking the
$0$-degree part of the Koszul resolution, provides a resolution of
$\mathcal{O}_{Y_g}$ on $\mathbb{P}(N_{X/Y}\oplus\mathcal{O}_Y)_g$.
\end{prop}

\section{Equivariant singular Bott-Chern classes}
Assume that $X$ is a $\mu_n$-equivariant projective manifold and
$S$ is a closed conical subset of $T_\R^*X_0$, we fix the
following notations:
\begin{displaymath}
\widetilde{\mathcal{U}}(X)=\bigoplus_{p\geq0}(D^{p,p}(X)/({\rm
Im}\partial+{\rm Im}\overline{\partial}))
\end{displaymath}
\begin{displaymath}
\widetilde{\mathcal{U}}(X,S)=\bigoplus_{p\geq0}(D^{p,p}(X,S)/({\rm
Im}\partial+{\rm Im}\overline{\partial})).
\end{displaymath}

\begin{defn}\label{501}
Let $i: Y\hookrightarrow X$ be an equivariant closed immersion of
projective manifolds. Let $N$ be the normal bundle of this
immersion and let $h_N$ be an invariant hermitian metric on $N$,
we shall denote $\overline{N}=(N,h_N)$. Moreover, let
$\overline{\eta}=(\eta,h_\eta)$ be an equivariant hermitian vector
bundle on $Y$ and let $\overline{\xi}.$ be a complex of
equivariant hermitian vector bundles on $X$ which provides a
resolution of $i_*\overline{\eta}$. The four-tuple
\begin{displaymath}
\overline{\Xi}=(i,\overline{N},\overline{\eta},\overline{\xi}.)
\end{displaymath} is called an equivariant hermitian embedded
vector bundle. Notice that an exact sequence of equivariant
hermitian vector bundles on $X$ is a particular case of
equivariant hermitian embedded vector bundle.
\end{defn}

\begin{defn}\label{502}
An equivariant singular Bott-Chern class for an equivariant
hermitian embedded vector bundle
$\overline{\Xi}=(i,\overline{N},\overline{\eta},\overline{\xi}.)$
is a class $\widetilde{H}\in \widetilde{\mathcal{U}}(X_g)$ such
that
\begin{displaymath}
{\rm dd}^c\widetilde{H}=\sum_{j}(-1)^j[{\rm
ch}_g(\overline{\xi}_j)]-{i_g}_*([{\rm ch}_g(\overline{\eta}){\rm
Td}_g^{-1}(\overline{N})]).
\end{displaymath}
\end{defn}

Note that the current
\begin{displaymath}
\sum_{j}(-1)^j[{\rm ch}_g(\overline{\xi}_j)]-{i_g}_*([{\rm
ch}_g(\overline{\eta}){\rm
Td}_g^{-1}(\overline{N})])=\sum_{j}(-1)^j[{\rm
ch}_g(\overline{\xi}_j)]-[{\rm ch}_g(\overline{\eta}){\rm
Td}_g^{-1}(\overline{N})]\delta_{Y_g}
\end{displaymath} is an
element in $D^*(X_g,N^\vee_{g,0})$, we would like to control the
singularities of the Bott-Chern class so that they are contained
in the same wave front set and we may do the pull-backs of
currents in certain situations. Theorem~\ref{310} allows us to do
this.

\begin{prop}\label{503}
Let
$\overline{\Xi}=(i,\overline{N},\overline{\eta},\overline{\xi}.)$
be an equivariant hermitian embedded vector bundle, then any
equivariant singular Bott-Chern class for $\overline{\Xi}$ belongs
to $\widetilde{\mathcal{U}}(X_g,N^\vee_{g,0})$.
\end{prop}
\begin{proof}
Firstly, note that Theorem~\ref{309} (ii) implies that the natural
map from $\widetilde{\mathcal{U}}(X_g,N^\vee_{g,0})$ to
$\widetilde{\mathcal{U}}(X_g)$ is injective. Then the statement in
this proposition does make sense and it follows from
Theorem~\ref{309} (i).
\end{proof}

Now assume that $f: X'\rightarrow X$ is an equivariant morphism of
projective manifolds which is transversal to $Y$. We formulate the
following Cartesian square
\begin{displaymath}
\xymatrix{Y' \ar[r]^-{i'} \ar[d]^h & X' \ar[d]^f \\
Y \ar[r]^-i & X.}
\end{displaymath} Since $h^*N$ is isomorphic to
the normal bundle of the immersion $i'$ (which implies that their
restrictions to the fixed point submanifolds are also isomorphic
to each other) and $f^*\xi.$ provides a resolution of
$i'_*h^*\eta$ on $X'$, we know that the notation
$f^*\overline{\Xi}=(i',h^*\overline{N},h^*\overline{\eta},f^*\overline{\xi}.)$
does make sense. Moreover, we conclude that $h_g^*N_g$ is
isomorphic to the normal bundle of $i'_g$. Then by
Proposition~\ref{503} and Theorem~\ref{304}, for any equivariant
singular Bott-Chern class $\widetilde{H}$ for $\overline{\Xi}$,
the pull-back $f_g^*\widetilde{H}$ is well-defined.

To every equivariant hermitian embedded vector bundle
$\overline{\Xi}=(i: Y\rightarrow
X,\overline{N},\overline{\eta},\overline{\xi}.)$, we may associate
two new equivariant hermitian embedded vector bundles. One is
${\rm tr}_1(\overline{\Xi}):=(j: \mathbb{P}_Y^1\rightarrow
W(i),\overline{N}_{W(i)/{\mathbb{P}_Y^1}},p_Y^*\overline{\eta},{\rm
tr}_1(\overline{\xi}.))$ concerning the construction of the
deformation to the normal cone, the other one is
$\overline{\Xi}_{Kos}:=(i_\infty: Y\rightarrow
\mathbb{P}(N\oplus\mathcal{O}_Y),\overline{N},\overline{\eta},K(\overline{\eta},\overline{N}))$
concerning the construction of the Koszul resolution.

Moreover, the direct sum of an equivariant hermitian embedded
vector bundle $\overline{\Xi}=(i: Y\rightarrow
X,\overline{N},\overline{\eta},\overline{\xi}.)$ with an exact
sequence $\overline{\varepsilon}$ of equivariant hermitian vector
bundles on $X$ is defined as
$\overline{\Xi}\oplus\overline{\varepsilon}:=(i: Y\rightarrow
X,\overline{N},\overline{\eta},\overline{\xi}.\oplus\overline{\varepsilon})$.

\begin{defn}\label{504}
Let $\Sigma$ be a set of equivariant hermitian embedded vector
bundles. We say that $\Sigma$ satisfies the condition (Hui) if

(i). any exact sequence of equivariant hermitian vector bundles on
an equivariant projective manifold belongs to $\Sigma$ and
$\Sigma$ is closed under the operation of taking direct sum with
an exact sequence of equivariant hermitian vector bundles;

(ii). for any element $\overline{\Xi}=(i: Y\rightarrow
X,\overline{N},\overline{\eta},\overline{\xi}.)\in \Sigma$ and for
every equivariant morphism $f: X'\rightarrow X$ of projective
manifolds which is transversal to $Y$, we have
$f^*\overline{\Xi}\in \Sigma$.

(iii). for any element $\overline{\Xi}=(i: Y\rightarrow
X,\overline{N},\overline{\eta},\overline{\xi}.)\in \Sigma$, the
associated equivariant hermitian embedded vector bundles ${\rm
tr}_1(\overline{\Xi})$ and $\overline{\Xi}_{Kos}$ both belong to
$\Sigma$.
\end{defn}

\begin{defn}\label{505}
Let $\Sigma$ be a set of equivariant hermitian embedded vector
bundles which satisfies the condition (Hui). A theory of
equivariant singular Bott-Chern classes for $\Sigma$ is an
assignment which, to each $\overline{\Xi}=(i: Y\rightarrow
X,\overline{N},\overline{\eta},\overline{\xi}.)\in \Sigma$,
assigns a class of currents
\begin{displaymath}
T(\overline{\Xi})\in \widetilde{\mathcal{U}}(X_g)
\end{displaymath} satisfying the following properties.

(i). (Differential equation) The following equality holds
\begin{displaymath}
{\rm dd}^cT(\overline{\Xi})=\sum_{j}(-1)^j[{\rm
ch}_g(\overline{\xi}_j)]-{i_g}_*([{\rm ch}_g(\overline{\eta}){\rm
Td}_g^{-1}(\overline{N})]).
\end{displaymath}

(ii). (Functoriality) For every equivariant morphism $f:
X'\rightarrow X$ of projective manifolds which is transversal to
$Y$, we have
\begin{displaymath}
f_g^*T(\overline{\Xi})=T(f^*\overline{\Xi}).
\end{displaymath}

(iii). (Normalization) Let $\overline{A}.$ be an equivariantly and
orthogonally split exact sequence of equivariant hermitian vector
bundles. Then
$T(\overline{\Xi})=T(\overline{\Xi}\oplus\overline{A}.)$.
Moreover, if $X={\rm Spec}(\C)$ is one point, $Y=\emptyset$ and
$\overline{\xi}.=0$, then $T(\overline{\Xi})=0$.
\end{defn}

\begin{rem}\label{506}
(i). When $Y=\emptyset$ and $\overline{\xi}.$ is an exact sequence
of equivariant hermitian vector bundles on $X$, the three
properties in the definition above imply that
\begin{displaymath}
T(\overline{\Xi})=\widetilde{{\rm ch}}_g(\overline{\xi}.)
\end{displaymath} where $\widetilde{{\rm ch}}_g$ is the equivariant
Bott-Chern secondary characteristic class associated to ${\rm
ch}_g$.

(ii). According to Definition~\ref{504}, the properties (ii) and
(iii) described in the definition above are reasonable.
\end{rem}

Throughout the rest of this section we shall assume that $\Sigma$
is a suitable set (big enough) of equivariant hermitian embedded
vector bundles and we shall also assume the existence of a theory
of equivariant singular Bott-Chern classes for $\Sigma$. We first
show the compatibility of equivariant singular Bott-Chern classes
with exact sequences and equivariant Bott-Chern secondary
characteristic classes.

We fix an equivariant closed immersion $i: Y\hookrightarrow X$ of
projective manifolds. Let
\begin{displaymath}
\overline{\chi}:\quad 0\rightarrow
\overline{\eta}_n\rightarrow\cdots\rightarrow\overline{\eta}_1\rightarrow\overline{\eta}_0\rightarrow0
\end{displaymath} be an exact sequence of equivariant hermitian
vector bundles on $Y$, and assume that we are given a family of
equivariant hermitian embedded vector bundles
$\{\overline{\Xi}_j=(i,\overline{N},\overline{\eta}_j,\overline{\xi}_{j,\cdot})\}_{j=0}^n$
which fit the following commutative diagram
\begin{displaymath}
\xymatrix{0 \ar[r] &
\overline{\xi}_{n,\cdot} \ar[r] \ar[d] & \cdots \ar[r] &
\overline{\xi}_{1,\cdot} \ar[r] \ar[d] &
\overline{\xi}_{0,\cdot} \ar[r] \ar[d] & 0 \\
0 \ar[r] & i_*\overline{\eta}_n \ar[r] & \cdots \ar[r] &
i_*\overline{\eta}_1 \ar[r] & i_*\overline{\eta}_0 \ar[r] & 0}
\end{displaymath} with exact rows. For each $k$, we write
$\overline{\varepsilon}_k$ for the exact sequence
\begin{displaymath}
0\rightarrow \overline{\xi}_{n,k}\rightarrow \cdots\rightarrow
\overline{\xi}_{1,k}\rightarrow \overline{\xi}_{0,k}\rightarrow 0.
\end{displaymath}

\begin{prop}\label{507}
Let notations and assumptions be as above, then we have the
following equality in $\widetilde{\mathcal{U}}(X_g)$
\begin{displaymath}
T(\bigoplus_{j\ even}\overline{\Xi}_j)-T(\bigoplus_{j\
odd}\overline{\Xi}_j)=\sum_k(-1)^k[\widetilde{{\rm
ch}}_g(\overline{\varepsilon}_k)]-{i_g}_*([{\rm
Td}_g^{-1}(\overline{N})\widetilde{{\rm ch}}_g(\overline{\chi})]).
\end{displaymath}
\end{prop}
\begin{proof}
According to Theorem~\ref{main} (v), we have the first
transgression exact sequences ${\rm tr}_1(\overline{\chi})$ on
$\mathbb{P}_Y^1$ and ${\rm tr}_1(\overline{\varepsilon}_k)$ on
$\mathbb{P}_X^1$ for each $k$. Denote by $l:
\mathbb{P}_Y^1\rightarrow \mathbb{P}_X^1$ the induced morphism,
then there exists an exact sequence of exact sequences
\begin{displaymath}
\cdots\rightarrow {\rm tr}_1(\overline{\varepsilon}_1)\rightarrow
{\rm tr}_1(\overline{\varepsilon}_0)\rightarrow l_*{\rm
tr}_1(\overline{\chi})\rightarrow 0.
\end{displaymath} We fix the
following notations
\begin{displaymath}
{\rm tr}_1(\overline{\chi})_+=\bigoplus_{j\ even}{\rm
tr}_1(\overline{\chi})_j,\quad\quad{\rm
tr}_1(\overline{\chi})_-=\bigoplus_{j\ odd}{\rm
tr}_1(\overline{\chi})_j, \end{displaymath}
\begin{displaymath}
{\rm tr}_1(\overline{\varepsilon}_k)_+=\bigoplus_{j\ even}{\rm
tr}_1(\overline{\varepsilon}_k)_j,\quad\quad{\rm
tr}_1(\overline{\varepsilon}_k)_-=\bigoplus_{j\ odd}{\rm
tr}_1(\overline{\varepsilon}_k)_j, \end{displaymath} then
\begin{displaymath}
{\rm tr}_1(\overline{\Xi})_+:=(l:
\mathbb{P}^1_Y\rightarrow\mathbb{P}^1_X,p_Y^*\overline{N},{\rm
tr}_1(\overline{\chi})_+,{\rm tr}_1(\overline{\varepsilon}.)_+),
\end{displaymath}
\begin{displaymath}
{\rm tr}_1(\overline{\Xi})_-:=(l:
\mathbb{P}^1_Y\rightarrow\mathbb{P}^1_X,p_Y^*\overline{N},{\rm
tr}_1(\overline{\chi})_-,{\rm tr}_1(\overline{\varepsilon}.)_-)
\end{displaymath} are two equivariant hermitian embedded vector
bundles.

By the functoriality of the first transgression exact sequences,
we obtain that
\begin{displaymath}
{\rm tr}_1(\overline{\Xi})_+\mid_{X\times\{0\}}=\bigoplus_{j\
even}{\rm tr}_1(\overline{\Xi}_j),\quad\quad{\rm
tr}_1(\overline{\Xi})_-\mid_{X\times\{0\}}=\bigoplus_{j\ odd}{\rm
tr}_1(\overline{\Xi}_j). \end{displaymath} Note that for any exact
sequence of equivariant hermitian vector bundles, its first
transgression exact sequence is equivariantly and orthogonally
split at infinity. So we have an isometry
\begin{displaymath}
{\rm tr}_1(\overline{\Xi})_+\mid_{X\times\{\infty\}}\cong{\rm
tr}_1(\overline{\Xi})_-\mid_{X\times\{\infty\}}.
\end{displaymath}
Since the wave front sets of the currents $[\log\mid z\mid^2]$ and
$T({\rm tr}_1(\overline{\Xi}_\pm))$ do not intersect (cf. \cite[P.
266]{BGS2}), by \cite[Thm. 8.2.10]{Hoe}, their products are
well-defined currents. Then in
$\widetilde{\mathcal{U}}(\mathbb{P}_{X_g}^1)$, we have
\begin{align*}
0&=\frac{\overline{\partial}}{2\pi i}\{\partial\log\mid
z\mid^2\cdot (T({\rm tr}_1(\overline{\Xi})_+)-T({\rm
tr}_1(\overline{\Xi})_-))\}+\frac{\partial}{2\pi i}\{\log\mid
z\mid^2\cdot\overline{\partial}(T({\rm
tr}_1(\overline{\Xi})_+)-T({\rm tr}_1(\overline{\Xi})_-))\}\\
&=(\frac{\overline{\partial}\partial}{2\pi i}\log\mid
z\mid^2)\cdot (T({\rm tr}_1(\overline{\Xi})_+)-T({\rm
tr}_1(\overline{\Xi})_-))-\log\mid
z\mid^2\cdot\frac{\overline{\partial}\partial}{2\pi i}(T({\rm
tr}_1(\overline{\Xi})_+)-T({\rm tr}_1(\overline{\Xi})_-))\\
&=(\delta_0-\delta_\infty)\cdot(T({\rm
tr}_1(\overline{\Xi})_+)-T({\rm tr}_1(\overline{\Xi})_-))-\log\mid
z\mid^2\cdot\sum_k(-1)^k({\rm ch}_g({\rm
tr}_1(\overline{\varepsilon}_k)_+)-{\rm ch}_g({\rm
tr}_1(\overline{\varepsilon}_k)_-))\\
&\qquad\qquad\qquad\qquad\qquad+\log\mid z\mid^2\cdot
{l_g}_*\{{\rm ch}_g({\rm tr}_1(\overline{\chi})_+){\rm
Td}_g^{-1}(p_Y^*\overline{N})-{\rm ch}_g({\rm
tr}_1(\overline{\chi})_-){\rm Td}_g^{-1}(p_Y^*\overline{N})\}.
\end{align*} Finally, integrating
both two sides of the equality above over $\mathbb{P}^1$ and using
the first construction of equivariant Bott-Chern secondary
classes, we get the identity in this proposition.
\end{proof}

A totally similar argument gives a proof of the following
proposition.

\begin{prop}\label{508}
Let
$\overline{\Xi}_0=(i,\overline{N}_0,\overline{\eta},\overline{\xi}.)$
be an equivariant hermitian embedded vector bundle with
$\overline{N}_0=(N,h_0)$. Assume that $h_1$ is another invariant
metric on $N$, we write $\overline{N}_1=(N,h_1)$ and
$\overline{\Xi}_1=(i,\overline{N}_1,\overline{\eta},\overline{\xi}.)$,
then we have
\begin{displaymath}
T(\overline{\Xi}_0)-T(\overline{\Xi}_1)=-{i_g}_*[{\rm
ch}_g(\overline{\eta})\widetilde{{\rm Td}_g^{-1}}(N,h_0,h_1)].
\end{displaymath}
\end{prop}

We now turn to a special case of closed immersion of equivariant
projective manifolds, namely the zero section embedding discussed
before Definition~\ref{402}. Precisely speaking, let $Y$ be an
equivariant projective manifold and let $\overline{\eta},
\overline{N}$ be two equivariant hermitian vector bundles on $Y$,
we denote $P=\mathbb{P}(N\oplus\mathcal{O}_Y)$. Let $\pi_P:
P\rightarrow Y$ be the canonical projection and let $i_\infty:
Y\rightarrow P$ be the zero section embedding. As in
Definition~\ref{402}, we shall write
$K(\overline{\eta},\overline{N})$ for the hermitian Koszul
resolution. We have already know that
$\overline{\Xi}_{Kos}(\overline{\eta},\overline{N})=(i_\infty,\overline{N},\overline{\eta},K(\overline{\eta},\overline{N}))$
is an equivariant hermitian embedded vector bundle associated to
$\overline{\Xi}$. Sometimes we just write it as
$K(\overline{\eta},\overline{N})$ for simplicity.

\begin{thm}\label{509}
Let $\Sigma$ be a set of equivariant hermitian embedded vector
bundles which satisfies the condition (Hui). Assume that $T$ is a
theory of equivariant singular Bott-Chern classes for $\Sigma$.
Then the current $({\pi_P}_g)_*T(K(\overline{\eta},\overline{N}))$
is ${\rm dd}^c$-closed. Moreover, the cohomology class that it
represents does not depend on the metrics on $\eta$ and $N$ so
that it determines a characteristic class $C_T(\eta,N)\in
\bigoplus_{p\geq0}H^{p,p}(Y_g)$.
\end{thm}
\begin{proof}
First note that the push-forwards for currents commute with
differentials by definition. Then we have
\begin{align*}
{\rm
dd}^c(({\pi_P}_g)_*T(K(\overline{\eta},\overline{N})))&=({\pi_P}_g)_*({\rm
dd^c}T(K(\overline{\eta},\overline{N})))\\
&=({\pi_P}_g)_*(\sum_k(-1)^k[{\rm
ch}_g(\wedge^k\overline{Q}^\vee)({\pi_P}_g)^*{\rm
ch}_g(\overline{\eta})]-{i_{\infty,g}}_*[{\rm
ch}_g(\overline{\eta}){\rm
Td}_g^{-1}(\overline{N})])\\
&=(({\pi_P}_g)_*[c_{{\rm rk}{Q_g}}(\overline{Q}_g){\rm
Td}_g^{-1}(\overline{Q})]-[{\rm
Td}_g^{-1}(\overline{N})])\cdot[{\rm ch}_g(\overline{\eta})].
\end{align*} We claim that $({\pi_P}_g)_*[c_{{\rm
rk}{Q_g}}(\overline{Q}_g){\rm Td}_g^{-1}(\overline{Q})]=[{\rm
Td}_g^{-1}(\overline{N})]$ so that
$({\pi_P}_g)_*T(K(\overline{\eta},\overline{N}))$ is ${\rm
dd}^c$-closed. Actually, one first need to notice that we have the
following tautological exact sequence on $P$
\begin{displaymath}
0\rightarrow \overline{\mathcal{O}(-1)}\rightarrow
\pi_P^*(\overline{N}\oplus\overline{\mathcal{O}_Y})\rightarrow\overline{Q}\rightarrow0.
\end{displaymath} Then, by restricting to the submanifold
$P_0=\mathbb{P}(N_g\oplus\mathcal{O}_{Y_g})$, we get a new exact
sequence
\begin{displaymath}
0\rightarrow \overline{\mathcal{O}(-1)}\mid_{P_0}\rightarrow
{\pi_P}_g^*(\overline{N}\mid_{Y_g}\oplus\overline{\mathcal{O}_{Y_g}})\mid_{P_0}\rightarrow\overline{Q}\mid_{P_0}\rightarrow0.
\end{displaymath} The $0$-degree part of the exact sequence above
is the tautological exact sequence on $P_0$. Taking the non-zero
degree part of this exact sequence we get an isometry
$\overline{Q}_\bot\mid_{P_0}\cong
(\pi_{P_0})^*(\overline{N}_\bot)$. Notice that the hermitian
complex $\wedge^\bullet \overline{Q}^\vee$ is equivariantly and
orthogonally split over $\mathbb{P}(N)$, so the support of
$c_{{\rm rk}{Q_g}}(\overline{Q}_g){\rm Td}_g^{-1}(\overline{Q})$
is contained in $P_0$. Moreover, by \cite[Cor. 3.8]{BL} we know
that
\begin{displaymath}
(\pi_{P_0})_*(c_{{\rm rk}{Q_g}}(\overline{Q}_g\mid_{P_0}){\rm
Td}^{-1}(\overline{Q}_g\mid_{P_0}))={\rm Td}^{-1}(\overline{N}_g).
\end{displaymath} Then we may compute
\begin{align*}
({\pi_P}_g)_*(c_{{\rm rk}{Q_g}}(\overline{Q}_g){\rm
Td}_g^{-1}(\overline{Q}))&=(\pi_{P_0})_*(c_{{\rm
rk}{Q_g}}(\overline{Q}_g){\rm
Td}_g^{-1}(\overline{Q}))\\
&=(\pi_{P_0})_*(c_{{\rm rk}{Q_g}}(\overline{Q}_g\mid_{P_0}){\rm
Td}^{-1}(\overline{Q}_g\mid_{P_0}){\rm
Td}_g^{-1}(\overline{Q}_\bot\mid_{P_0}))\\
&=(\pi_{P_0})_*(c_{{\rm rk}{Q_g}}(\overline{Q}_g\mid_{P_0}){\rm
Td}^{-1}(\overline{Q}_g\mid_{P_0}){\rm
Td}_g^{-1}((\pi_{P_0})^*(\overline{N}_\bot))\\
&={\rm Td}^{-1}(\overline{N}_g){\rm
Td}_g^{-1}(\overline{N}_\bot)={\rm Td}_g^{-1}(\overline{N})
\end{align*} which completes the proof of the claim.

Now let $h_0$ and $h_1$ (resp. $g_0$ and $g_1$) be two invariant
hermitian metrics on $N$ (resp. $\eta$). We write
$\overline{N}_i=(N,h_i)$ and $\overline{\eta}_i=(\eta,g_i)$. We
denote also by $h_0$ and $h_1$ the metrics induced on $Q^\vee$.
Then by Proposition~\ref{507}, we have
\begin{align*}
&({\pi_P}_g)_*(T(K(\overline{\eta}_0,\overline{N}_0))-T(K(\overline{\eta}_1,\overline{N}_0)))\\
=&({\pi_P}_g)_*[\sum_k(-1)^k{\rm
ch}_g(\wedge^k\overline{Q}_0^\vee){\pi_P}_g^*\widetilde{{\rm
ch}}_g(\eta,g_0,g_1)]-[{\rm
Td}_g^{-1}(\overline{N}_0)\widetilde{{\rm ch}}_g(\eta,g_0,g_1)].
\end{align*}
So using the projection formula and our claim before, we get
\begin{displaymath}
({\pi_P}_g)_*(T(K(\overline{\eta}_0,\overline{N}_0))=({\pi_P}_g)_*(T(K(\overline{\eta}_1,\overline{N}_0)).
\end{displaymath}

On the other hand, applying Proposition~\ref{507} and
Proposition~\ref{508}, we have
\begin{align*}
&({\pi_P}_g)_*(T(K(\overline{\eta}_1,\overline{N}_0))-T(K(\overline{\eta}_1,\overline{N}_1)))\\
=&({\pi_P}_g)_*[\sum_k(-1)^k\widetilde{{\rm
ch}}_g(\wedge^kQ^\vee,h_0,h_1){\pi_P}_g^*{\rm
ch}_g(\overline{\eta}_1)]-[{\rm
ch}_g(\overline{\eta}_1)\widetilde{{\rm
Td}_g^{-1}}(N,h_0,h_1)]\\
=&\{({\pi_P}_g)_*\sum_k(-1)^k[\widetilde{{\rm
ch}}_g(\wedge^kQ^\vee,h_0,h_1)]-[\widetilde{{\rm
Td}_g^{-1}}(N,h_0,h_1)]\}\cdot[{\rm ch}_g(\overline{\eta}_1)].
\end{align*} We construct the first transgression exact sequence of
$0\rightarrow 0\rightarrow (N,h_1)\rightarrow (N,h_0)\rightarrow
0$ on $\mathbb{P}_Y^1$ so that we may have an equivariant
hermitian vector bundle $(\widetilde{N},h^{\widetilde{N}})$ on
$\mathbb{P}_Y^1$ such that
\begin{displaymath}
(\widetilde{N},h^{\widetilde{N}})\mid_{Y\times\{0\}}=(N,h_1),\quad\quad
(\widetilde{N},h^{\widetilde{N}})\mid_{Y\times\{\infty\}}=(N,h_0).
\end{displaymath} Now we apply the Koszul construction to the
bundles $p_Y^*\overline{\eta}_1$ and
$(\widetilde{N},h^{\widetilde{N}})$ and denote by $\pi_W$ the
canonical projection from
$W:=\mathbb{P}(\widetilde{N}\oplus\mathcal{O}_{\mathbb{P}_Y^1})$
to $\mathbb{P}_Y^1$. By the universal properties of projective
space bundle and fibre product,
$\mathbb{P}(\widetilde{N}\oplus\mathcal{O}_{\mathbb{P}_Y^1})=\mathbb{P}(p_Y^*N\oplus
p_Y^*\mathcal{O}_Y)$ which is isomorphic to
$\mathbb{P}(N\oplus\mathcal{O}_Y)\times \mathbb{P}^1$ and the
tautological quotient bundle on $W$ is isomorphic to
$\widetilde{Q}$ whose definition is similar to that of
$\widetilde{N}$. Thus we have
\begin{align*}
&({\pi_P}_g)_*\sum_k(-1)^k[\widetilde{{\rm
ch}}_g(\wedge^kQ^\vee,h_0,h_1)]\\
=&({\pi_P}_g)_*\sum_k(-1)^k({p_P}_g)_*([-\log\mid
z\mid^2]\cdot[{\rm
ch}_g(\wedge^k(\widetilde{Q},h^{\widetilde{Q}}))])\\
=&({p_Y}_g)_*({\pi_W}_g)_*\sum_k(-1)^k([-\log\mid
z\mid^2]\cdot[{\rm
ch}_g(\wedge^k(\widetilde{Q},h^{\widetilde{Q}}))])\\
=&({p_Y}_g)_*([-\log\mid z\mid^2]\cdot[{\rm
Td}_g^{-1}(\widetilde{N},h^{\widetilde{N}})])=[\widetilde{{\rm
Td}_g^{-1}}(N,h_0,h_1)].
\end{align*} So we have proved that
$({\pi_P}_g)_*T(K(\overline{\eta},\overline{N}))$ dose not depend
on the choices of the metrics. Thus we have a well-defined class
$C_T(\eta,N)$. The fact that this characteristic class
$C_T(\eta,N)$ belongs to $\bigoplus_{p\geq0}H^{p,p}(Y_g)$ follows
from \cite[Theroem 1.2.2 (iii)]{GS2}.
\end{proof}

\section{Classification of theories of equivariant singular
Bott-Chern classes} The aim of this section is to give some
results concerning the classification of all possible theories of
equivariant singular Bott-Chern classes. We shall prove that a
theory of equivariant singular Bott-Chern classes $T$ is totally
determined by the characteristic class $C_T$ introduced in last
section. Our main theorem is the following.

\begin{thm}\label{601}
Let $\Sigma$ be a set of equivariant hermitian embedded vector
bundles which satisfies the condition (Hui). Assume that $T$ and
$T'$ are two theories of equivariant singular Bott-Chern classes
for $\Sigma$. Then $T=T'$ if and only if for any
$(i,\overline{N},\overline{\eta},\overline{\xi}.)\in \Sigma$,
$C_T(\eta,N)=C_{T'}(\eta,N)$.
\end{thm}
\begin{proof}
One direction is clear. For the other one, we assume that
$C_T=C_{T'}$. Let $\overline{\Xi}=(i: Y\rightarrow
X,\overline{N},\overline{\eta},\overline{\xi}.)$ be an element in
$\Sigma$. As before, we denote by $W$ the deformation to the
normal cone and denote by $p_W$ the composition of $p_X$ and the
blow-down map $\pi$. Moreover, we write $p_{\widetilde{X}}:
\widetilde{X}\rightarrow X$ and $p_P: P\rightarrow X$ for the
morphisms induced by $p_W$. The morphism $p_P$ can be factored as
$i\circ \pi_P$.

The normal bundle of the immersion $j: Y\times
\mathbb{P}^1\rightarrow W$ is isomorphic to $p_Y^*N\otimes
q_Y^*\mathcal{O}(-1)$. We endow it with the hermitian metric
induced by the metric on $N$ and the Fubini-Study metric on
$\mathcal{O}(-1)$, the corresponding hermitian vector bundle will
be denoted by $\overline{N}'$.

By Theorem~\ref{main}, the restriction of ${\rm
tr}_1(\overline{\xi}.)$ to $X\times\{0\}$ is isometric to
$\overline{\xi}.$ and the restriction of ${\rm
tr}_1(\overline{\xi}.)$ to $\widetilde{X}$ is equivariantly and
orthogonally split. Moreover, the restriction of ${\rm
tr}_1(\overline{\xi}.)$ to $P$ fits an exact sequence
\begin{displaymath}
0\rightarrow \overline{A}.\rightarrow {\rm
tr}_1(\overline{\xi}.)\mid_P\rightarrow
K(\overline{\eta},\overline{N})\rightarrow 0
\end{displaymath}
where $\overline{A}.$ is an equivariantly and orthogonally split
exact sequence. We denote by $\overline{\varepsilon}_k$ the exact
sequence of the following exact sequence of equivariant hermitian
vector bundles
\begin{displaymath}
0\rightarrow \overline{A}_k\rightarrow {\rm
tr}_1(\overline{\xi}.)_k\mid_P\rightarrow
K(\overline{\eta},\overline{N})_k\rightarrow 0.
\end{displaymath}
Next, we write $U$ for the current $[-\log\mid z\mid^2]$ on
$\mathbb{P}^1$ associated to a locally integrable differential
form. Its pull-back to $W_g$ is also locally integrable hence
defines a current on $W_g$ which will be also denoted by $U$. Note
that $q_{W(i_g)}=q_{W_g}\circ i_{W(i_g)}$ where $i_{W(i_g)}$ is
the natural immersion from $W(i_g)$ to $W_g$ and the wave front
set of $T({\rm tr}_1(\overline{\Xi}))$ is contained in the
conormal bundle ${N'}_g^\vee$. Hence the wave front sets of $U$
and $T({\rm tr}_1(\overline{\Xi}))$ are disjoint so that their
product $U\cdot T({\rm tr}_1(\overline{\Xi}))$ is a well-defined
current on $W_g$. Then, using the properties of equivariant
singular Bott-Chern classes in Definition~\ref{505}, the equality
\begin{align*}
0&={\rm dd}^c({p_W}_g)_*(U\cdot T({\rm
tr}_1(\overline{\Xi})))\\
&=({p_{\widetilde{X}}}_g)_*(T({\rm
tr}_1(\overline{\Xi}))\mid_{\widetilde{X}_g})+({p_P}_g)_*(T({\rm
tr}_1(\overline{\Xi}))\mid_{P_g})-T(\overline{\Xi})\\
&\qquad\qquad\qquad-({p_W}_g)_*(U\cdot(\sum_k(-1)^k{\rm ch}_g({\rm
tr}_1(\overline{\xi}.)_k)-{j_g}_*({\rm
ch}_g(p_Y^*\overline{\eta}){\rm Td}_g^{-1}(\overline{N}'))))
\end{align*}
holds in $\widetilde{\mathcal{U}}(X_g)$. Notice that $T({\rm
tr}_1(\overline{\Xi}))\mid_{\widetilde{X}_g}=T({\rm
tr}_1(\overline{\Xi})\mid_{\widetilde{X}})=\widetilde{{\rm
ch}}_g({\rm tr}_1(\overline{\xi}.)\mid_{\widetilde{X}})=0$ and by
Proposition~\ref{507}, we have
\begin{displaymath}
T({\rm tr}_1(\overline{\Xi}))\mid_{P_g}=T({\rm
tr}_1(\overline{\Xi})\mid_P)=T(K(\overline{\eta},\overline{N}))-\sum_k(-1)^k[\widetilde{{\rm
ch}}_g(\overline{\varepsilon}_k)].
\end{displaymath} Moreover,
using the factorization of $p_P$, we have
\begin{displaymath}
({p_P}_g)_*T(K(\overline{\eta},\overline{N}))={i_g}_*({\pi_P}_g)_*T(K(\overline{\eta},\overline{N}))={i_g}_*C_T(\eta,N).
\end{displaymath} By the properties of the Fubini-Study metric,
${\rm ch}_g(p_Y^*\overline{\eta}){\rm Td}_g^{-1}(\overline{N}')$
is invariant under the involution on $\mathbb{P}^1$ which sends
$z$ to $1/z$. Thus we obtain
\begin{displaymath}
({p_W}_g)_*(U\cdot ({j_g}_*({\rm ch}_g(p_Y^*\overline{\eta}){\rm
Td}_g^{-1}(\overline{N}'))))={i_g}_*({p_Y}_g)_*(U\cdot ({\rm
ch}_g(p_Y^*\overline{\eta}){\rm Td}_g^{-1}(\overline{N}')))=0
\end{displaymath}
since the current $U$ really changes its sign under the involution
$z\rightarrow 1/z$. Gathering all computations above we finally
get the following current equation
\begin{displaymath}
T(\overline{\Xi})=-({p_W}_g)_*(U\cdot\sum_k(-1)^k{\rm ch}_g({\rm
tr}_1(\overline{\xi}.)_k))-\sum_k(-1)^k({p_P}_g)_*[\widetilde{{\rm
ch}}_g(\overline{\varepsilon}_k)]+{i_g}_*C_T(\eta,N).
\end{displaymath} A similar current equation for $T'$ can be
obtained in the same way. By our assumption we have
$C_T(\eta,N)=C_{T'}(\eta,N)$, so that
$T(\overline{\Xi})=T'(\overline{\Xi})$. This completes the whole
proof.
\end{proof}

From the proof of Theorem~\ref{601}, it is natural to guess that
if we are given an explicit definition of the equivariant
characteristic class $C$, we then get a theory of equivariant
singular Bott-Chern classes $T$ such that $C_T$ is exactly $C$. We
prove this conjecture in the following theorem.

\begin{thm}\label{602}
Let $\Sigma$ be a set of equivariant hermitian embedded vector
bundles which satisfies the condition (Hui). Assume that $C$ is an
equivariant characteristic class for pairs of vector bundles
$(\eta,N)$ which appear in the elements
$(i,\overline{N},\overline{\eta},\overline{\xi}.)\in \Sigma$. Then
there exists a theory of equivariant singular Bott-Chern classes
for $\Sigma$ such that $C_T=C$.
\end{thm}
\begin{proof}
For any element $\overline{\Xi}=(i: Y\rightarrow
X,\overline{N},\overline{\eta},\overline{\xi}.)\in \Sigma$, we
define
\begin{displaymath}
T(\overline{\Xi})=-({p_W}_g)_*(U\cdot\sum_k(-1)^k{\rm ch}_g({\rm
tr}_1(\overline{\xi}.)_k))-\sum_k(-1)^k({p_P}_g)_*[\widetilde{{\rm
ch}}_g(\overline{\varepsilon}_k)]+{i_g}_*C(\eta,N).
\end{displaymath} Our first aim is to prove that such $T$ does not
depend on the choice of the metric on ${\rm tr}_1(\xi.)$ or on
$A.$ and that such $T$ satisfies all properties in the definition
of a theory of equivariant singular Bott-Chern classes.

We denote by $h_k$ and $h'_k$ (resp. $g_k$ and $g'_k$) two
invariant hermitian metrics on ${\rm tr}_1(\xi.)_k$ (resp. $A_k$)
such that the resulting hermitian vector bundles all satisfy the
requirements in Theorem~\ref{main}. Then, in
$\widetilde{\mathcal{U}}(X_g)$, we have
\begin{align*}
&\sum_k(-1)^k({p_P}_g)_*[\widetilde{{\rm
ch}}_g(\overline{\varepsilon}_k)]-\sum_k(-1)^k({p_P}_g)_*[\widetilde{{\rm
ch}}_g(\overline{\varepsilon}'_k)]\\
=&\sum_k(-1)^k({p_P}_g)_*[\widetilde{{\rm
ch}}_g(A_k,g_k,g'_k)]-\sum_k(-1)^k({p_P}_g)_*[\widetilde{{\rm
ch}}_g({\rm tr}_1(\xi.)_k\mid_P,h_k,h'_k)].
\end{align*} The first
term of the right-hand side vanishes due to Proposition~\ref{507}
and the assumption that the complex $A.$ is orthogonally split for
both metrics.

On the other hand, we have by definition
\begin{align*}
&({p_W}_g)_*(U\cdot\sum_k(-1)^k{\rm ch}_g({\rm
tr}_1(\xi.)_k,h_k))-({p_W}_g)_*(U\cdot\sum_k(-1)^k{\rm ch}_g({\rm
tr}_1(\xi.)_k,h'_k))\\
=&({p_W}_g)_*(U\cdot\sum_k(-1)^k{\rm dd}^c\widetilde{{\rm
ch}}_g({\rm tr}_1(\xi.)_k,h_k,h'_k)).
\end{align*} But, in
$\widetilde{\mathcal{U}}(X_g)$, we have
\begin{align*}
&({p_W}_g)_*(U\cdot\sum_k(-1)^k{\rm dd}^c\widetilde{{\rm
ch}}_g({\rm tr}_1(\xi.)_k,h_k,h'_k))\\
=&\sum_k(-1)^k({p_{\widetilde{X}}}_g)_*[\widetilde{{\rm
ch}}_g({\rm
tr}_1(\xi.)_k,h_k,h'_k)]\mid_{\widetilde{X}_g}+\sum_k(-1)^k({p_P}_g)_*[\widetilde{{\rm
ch}}_g({\rm tr}_1(\xi.)_k,h_k,h'_k)]\mid_{P_g}\\
&\qquad\qquad\qquad\qquad\qquad\qquad-\sum_k(-1)^k[\widetilde{{\rm
ch}}_g({\rm tr}_1(\xi.)_k,h_k,h'_k)]\mid_{X\times\{0\}}.
\end{align*} The last
term of the right-hand side vanishes because the metrics $h_k$ and
$h'_k$ agree each other on $X\times\{0\}$. The first term vanishes
due to the assumption that ${\rm tr}_1(\xi.)\mid_{\widetilde{X}}$
is orthogonally split with both metrics. Therefore, combining the
two computations above, we know that the definition of $T$ is
independent of the metrics on $A.$ and ${\rm tr}_1(\xi.)$.

We next prove that the definition of $T$ satisfies the three
properties in the definition of a theory of equivariant singular
Bott-Chern classes. For the differential equation, we compute
\begin{align*}
{\rm
dd}^cT(\overline{\Xi})=&-\sum_k(-1)^k({p_{\widetilde{X}}}_g)_*{\rm
ch}_g({\rm
tr}_1(\overline{\xi}.)_k\mid_{\widetilde{X}})-\sum_k(-1)^k({p_P}_g)_*{\rm
ch}_g({\rm tr}_1(\overline{\xi}.)_k\mid_P)\\
&+\sum_k(-1)^k{\rm ch}_g({\rm
tr}_1(\overline{\xi}.)_k\mid_{X\times\{0\}})\\
&-\sum_k(-1)^k({p_P}_g)_*({\rm ch}_g(\overline{A}_k)+{\rm
ch}_g(K(\overline{\eta},\overline{N})_k)-{\rm ch}_g({\rm
tr}_1(\overline{\xi}.)_k\mid_P)). \end{align*} Using the fact that
$\overline{A}.$ and ${\rm
tr}_1(\overline{\xi}.)\mid_{\widetilde{X}}$ are equivariantly and
orthogonally split we obtain
\begin{align*}
{\rm dd}^cT(\overline{\Xi})&=\sum_k(-1)^k{\rm
ch}_g(\overline{\xi}_k)-\sum_k(-1)^k({p_P}_g)_*{\rm
ch}_g(K(\overline{\eta},\overline{N})_k)\\
&=\sum_k(-1)^k[{\rm ch}_g(\overline{\xi}_k)]-({p_P}_g)_*[c_{{\rm
rk}Q_g}(\overline{Q}_g){\rm Td}_g^{-1}(\overline{Q}){\rm
ch}_g(\pi_P^*\overline{\eta})]\\
&=\sum_k(-1)^k[{\rm ch}_g(\overline{\xi}_k)]-{i_g}_*[{\rm
ch}_g(\overline{\eta}){\rm Td}_g^{-1}(\overline{N})].
\end{align*}

Secondly, the functoriality property for our definition of $T$
follows from the functoriality property for ${\rm ch}_g$,
$\widetilde{{\rm ch}}_g$ and $C$.

We now prove the normalization property. We first assume that
$Y=\emptyset$ and $\overline{\xi}.$ is an equivariantly and
orthogonally split exact sequence. This means that if we write
$\overline{K}_i={\rm Ker}({\rm d}_i: \overline{\xi}_i\rightarrow
\overline{\xi}_{i-1})$, then $\overline{\xi}_i$ is isometric to
$\overline{K}_i\oplus\overline{K}_{i-1}$. Hence by the
construction of ${\rm tr}_1(\xi.)$, we know that
\begin{displaymath}
{\rm tr}_1(\xi.)_i=p_X^*K_i\otimes q_X^*\mathcal{O}(i)\oplus
p_X^*K_{i-1}\otimes q_X^*\mathcal{O}(i-1).
\end{displaymath} This
formula implies that $\sum_k(-1)^k{\rm ch}_g({\rm
tr}_1(\overline{\xi}.)_k)$ is invariant under the involution on
$\mathbb{P}^1$ which sends $z$ to $1/z$. So the first term in the
definition for $T$ vanishes. It is clear that the other two terms
also vanish in this special case. Hence we obtain
$T(\overline{\xi}.)=0$. Now let $\overline{\Xi}=(i: Y\rightarrow
X,\overline{N},\overline{\eta},\overline{\xi}.)$ and let
$\overline{B}.$ be an equivariantly and orthogonally split exact
sequence of equivariant hermitian vector bundles on $X$. By
\cite[Section 1.1]{GS5}, we have
\begin{displaymath}
{\rm tr}_1(\xi.\oplus B.)={\rm tr}_1(\xi.)\oplus \pi^*{\rm
tr}_1(B.). \end{displaymath} In order to compute
$T(\overline{\Xi}\oplus \overline{B}.)$, we consider the following
exact sequences
\begin{displaymath}
\overline{\varepsilon}'_k:\quad
0\rightarrow\overline{A}_k\oplus\pi^*{\rm
tr}_1(\overline{B}.)_k\mid_P\rightarrow {\rm
tr}_1(\overline{\xi}.)_k\oplus \pi^*{\rm
tr}_1(\overline{B}.)_k\mid_P\rightarrow
K(\overline{\eta},\overline{N})_k\rightarrow 0. \end{displaymath}
By the additivity of equivariant Bott-Chern secondary
characteristic classes, we have $\widetilde{{\rm
ch}}_g(\overline{\varepsilon}_k)=\widetilde{{\rm
ch}}_g(\overline{\varepsilon}'_k)$. Again using the additivity of
equivariant Chern classes, we finally get
\begin{displaymath}
T(\overline{\Xi}\oplus\overline{B}.)-T(\overline{B}.)=0.
\end{displaymath}

At last, we should prove that the equivariant characteristic class
$C_T$ is exactly equal to $C$. Note that the arguments above show
that $T$ is really a theory of equivariant singular Bott-Chern
classes, then as what we have seen in the proof of
Theorem~\ref{601}, for any element $\overline{\Xi}=(i:
Y\rightarrow X,\overline{N},\overline{\eta},\overline{\xi}.)\in
\Sigma$ we always have
\begin{displaymath}
T(\overline{\Xi})=-({p_W}_g)_*(U\cdot\sum_k(-1)^k{\rm ch}_g({\rm
tr}_1(\overline{\xi}.)_k))-\sum_k(-1)^k({p_P}_g)_*[\widetilde{{\rm
ch}}_g(\overline{\varepsilon}_k)]+{i_g}_*C_T(\eta,N).
\end{displaymath} In particular, for the Koszul construction $(i:
Y\rightarrow
\mathbb{P}(N\oplus\mathcal{O}_Y),\overline{N},\overline{\eta},K(\overline{\eta},\overline{N}))$,
we have
\begin{displaymath}
T(K(\overline{\eta},\overline{N}))=-({p_W}_g)_*(U\cdot\sum_k(-1)^k{\rm
ch}_g({\rm
tr}_1(K(\overline{\eta},\overline{N})_k))-\sum_k(-1)^k({p_P}_g)_*[\widetilde{{\rm
ch}}_g(\overline{\varepsilon}_k)]+{i_g}_*C_T(\eta,N).
\end{displaymath} Comparing with the definition of $T$ via the
characteristic class $C$, we get
${i_g}_*C_T(\eta,N)={i_g}_*C(\eta,N)$ and hence
$C_T(\eta,N)=C(\eta,N)$ after composing $({\pi_P}_g)_*$. This
completes the whole proof.
\end{proof}

To end this section, we shall give an example of the set of
equivariant hermitian embedded vector bundles which satisfies the
condition (Hui) and we shall also give a general way to construct
the characteristic class $C$.

\begin{defn}\label{603}
Let $\overline{\Xi}=(i: Y\rightarrow
X,\overline{N},\overline{\eta},\overline{\xi.})$ be an equivariant
hermitian embedded vector bundle. The equivariant rank of
$\overline{\Xi}$ is the sequence of locally constant functions
$({\rm rk}\eta_\zeta)_{\zeta\in S^1}$. The equivariant codimension
of $\overline{\Xi}$ is the sequence of locally constant functions
$({\rm rk}N_\zeta)_{\zeta\in S^1}$. When $Y=\emptyset$, we shall
say that an exact sequence of equivariant hermitian vector bundles
on $X$ has arbitrary equivariant rank and arbitrary equivariant
codimension.
\end{defn}

\begin{prop}\label{604}
Let $(t_\zeta)_{\zeta\in S^1}$ and $(s_\zeta)_{\zeta\in S^1}$ be
two sequences of natural numbers. Let $\Sigma$ be a set consisting
of all equivariant hermitian embedded vector bundles of
equivariant rank less than or equal to $(t_\zeta)_{\zeta\in S^1}$
and of equivariant codimension less than or equal to
$(s_\zeta)_{\zeta\in S^1}$. Then $\Sigma$ satisfies the condition
(Hui).
\end{prop}
\begin{proof}
The first requirement in the condition (Hui) is naturally
fulfilled by definition. For the second requirement, let $(i:
Y\rightarrow X,\overline{N},\overline{\eta},\overline{\xi.})$ be
an element in $\Sigma$ and let $f: X'\rightarrow X$ be an
equivariant morphism which is transversal to $Y$. Then $f^{-1}(Y)$
either is empty set or has the same codimension as $Y$. In the
first case, we are done. In the second case, we formulate the
following Cartesian square
\begin{displaymath}
\xymatrix{ Y' \ar[r]^-{i'} \ar[d]^h & X' \ar[d]^f \\
Y \ar[r]^-{i} & X,} \end{displaymath} then $h^*N\cong N'$. Note
that $(h^*\eta)\mid_{Y'_g}=h_g^*(\eta\mid_{Y_g})$ and
$(h^*N)\mid_{Y'_g}=h_g^*(N\mid_{Y_g})$, we have the inequalities
$({\rm rk}(h^*\eta)_\zeta)_{\zeta\in S^1}\leq(t_\zeta)_{\zeta\in
S^1}$ and $({\rm rk}N'_\zeta)_{\zeta\in
S^1}\leq(s_\zeta)_{\zeta\in S^1}$. This means that the equivariant
hermitian embedded vector bundle $(i': Y'\rightarrow
X',\overline{N'},h^*\overline{\eta},f^*\overline{\xi.})$ is also
in $\Sigma$. For the last requirement in the condition (Hui), we
again let $\overline{\Xi}=(i: Y\rightarrow
X,\overline{N},\overline{\eta},\overline{\xi.})$ be an element in
$\Sigma$. Then the associated Koszul construction
$\overline{\Xi}_{Kos}$ clearly has the same equivariant rank and
equivariant codimension as $\overline{\Xi}$. Concerning the
construction of the deformation to the normal cone, ${\rm
tr}_1(\overline{\Xi})$ clearly has the same equivariant rank as
$\overline{\Xi}$. Moreover, the normal bundle of $Y\times
\mathbb{P}^1$ in $W(i)$ is $N'=p_Y^*N\otimes
q_Y^*\mathcal{O}(-1)$. Note that we assume that $\mathbb{P}^1$
admits the trivial $g$-action, then $({\rm rk}N_\zeta)_{\zeta\in
S^1}=({\rm rk}N'_\zeta)_{\zeta\in S^1}$ so that ${\rm
tr}_1(\overline{\Xi})$ also has the same equivariant codimension
as $\overline{\Xi}$. Therefore, we have that ${\rm
tr}_1(\overline{\Xi})$ and $\overline{\Xi}_{Kos}$ are both
elements in $\Sigma$.
\end{proof}

We finally give a general construction of the characteristic class
$C$ for the set $\Sigma$ in last proposition.

\begin{defn}\label{605}
Let $(\varphi_\zeta)_{\zeta\in S^1}$ be a family of
$\GL(\C)$-invariant formal power series such that
$\varphi_\zeta\in \C[[\mathbf{gl}_{{\rm rk}\eta_\zeta}(\C)]]$. And
let $(\psi_\zeta)_{\zeta\in S^1}$ be a family of
$\GL(\C)$-invariant formal power series such that $\psi_\zeta\in
\C[[\mathbf{gl}_{{\rm rk}N_\zeta}(\C)]]$. Moreover, let $\phi\in
\C[[\bigoplus_{\zeta\in S^1}\C\oplus\bigoplus_{\zeta\in S^1}\C]]$
be any formal power series. We define the equivariant character
form $\phi_g(\overline{\eta},\overline{N})$ as
\begin{displaymath}
\phi_g(\overline{\eta},\overline{N})=\phi((\varphi_{\zeta}(-\frac{\Omega^{\overline{\eta}_\zeta}}{2\pi
i}))_{\zeta\in
S^1},(\psi_{\zeta}(-\frac{\Omega^{\overline{N}_\zeta}}{2\pi
i}))_{\zeta\in S^1}). \end{displaymath} The cohomology class that
$\phi_g(\overline{\eta},\overline{N})$ represents is independent
of the choices of the metrics, hence it define a characteristic
class. We denote it by $C(\eta,N)$.
\end{defn}

Then the following corollary follows immediately from
Theorem~\ref{601} and Theorem~\ref{602}.

\begin{cor}\label{606}
Let $\Sigma$ be the set of equivariant hermitian embedded vector
bundles defined in Proposition~\ref{604}. Let $C$ be an
equivariant characteristic class for pairs of equivariant
hermitian vector bundles given in the way as in
Definition~\ref{605}. Then there exists a unique theory of
equivariant singular Bott-Chern classes $T$ for $\Sigma$ such that
$C_T$ is equal to $C$.
\end{cor}

\section{Compatibility with the projection formula}
As usual, let $\Sigma$ be a set of equivariant hermitian embedded
vector bundles which satisfies the condition (Hui). In this
section, we shall give the sufficient and necessary condition for
a theory of equivariant singular Bott-Chern classes to be
compatible with the projection formula. This can be regarded as an
example of how the properties of the characteristic class $C_T$
are reflected in the corresponding theory of equivariant singular
Bott-Chern classes.

Now, let $\overline{\Xi}=(i: Y\rightarrow
X,\overline{N},\overline{\eta},\overline{\xi}.)$ be an equivariant
hermitian embedded vector bundle in $\Sigma$. For any equivariant
hermitian vector bundle $\overline{\kappa}$ on $X$, we define
\begin{displaymath}
\overline{\Xi}\otimes \overline{\kappa}=(i: Y\rightarrow
X,\overline{N},\overline{\eta}\otimes
i^*\overline{\kappa},\overline{\xi}.\otimes \overline{\kappa}).
\end{displaymath}
Note that $\overline{\Xi}\otimes \overline{\kappa}$ is also an
equivariant hermitian embedded vector bundle according to the
projection formula. We assume that $\Sigma$ is big enough so that
$\overline{\Xi}\otimes \overline{\kappa}$ and all equivariant
hermitian embedded vector bundle appearing below belong to it.

\begin{defn}\label{701}
Let notations and assumptions be as above. A theory of equivariant
singular Bott-Chern classes $T$ for $\Sigma$ is said to be
compatible with the projection formula if
\begin{displaymath}
T(\overline{\Xi}\otimes \overline{\kappa})=T(\overline{\Xi})\cdot
{\rm ch}_g(\overline{\kappa}). \end{displaymath}
\end{defn}

\begin{prop}\label{702}
Let notations and assumptions be as above. Then
\begin{displaymath}
T(\overline{\Xi}\otimes \overline{\kappa})-T(\overline{\Xi})\cdot
{\rm ch}_g(\overline{\kappa})={i_g}_*(C_T(\eta\otimes
i^*\kappa,N))-{i_g}_*(C_T(\eta,N))\cdot{\rm
ch}_g(\overline{\kappa}). \end{displaymath}
\end{prop}
\begin{proof}
As before, denote by $p_W$ the composition of the blow-down map
$\pi$ and the projection $p_X: X\times \mathbb{P}^1\rightarrow X$.
Then by the construction of ${\rm tr}_1(\cdot)$, we have ${\rm
tr}_1(\overline{\xi}.\otimes \overline{\kappa})={\rm
tr}_1(\overline{\xi}.)\otimes p_W^*\overline{\kappa}$. Then, on
one hand, we have
\begin{align*}
({p_W}_g)_*(U\cdot\sum_k(-1)^k{\rm ch}_g({\rm
tr}_1(\overline{\xi}.\otimes
\overline{\kappa})_k))&=({p_W}_g)_*(U\cdot\sum_k(-1)^k{\rm
ch}_g({\rm tr}_1(\overline{\xi}.)_k){p_W}_g^*{\rm
ch}_g(\overline{\kappa}))\\
&=({p_W}_g)_*(U\cdot\sum_k(-1)^k{\rm ch}_g({\rm
tr}_1(\overline{\xi}.)_k)){\rm ch}_g(\overline{\kappa}).
\end{align*} On the other hand, the Koszul resolution of
$i_*(\eta\otimes i^*\kappa)$ is given by
\begin{displaymath}
K(\eta\otimes i^*\kappa,N)=K(\eta,N)\otimes p_P^*\kappa.
\end{displaymath} Then for each $k$, if we write
$\overline{\varepsilon}_k\otimes p_P^*\overline{\kappa}$ for the
exact sequence
\begin{displaymath}
0\rightarrow \overline{A}_k\otimes
p_P^*\overline{\kappa}\rightarrow {\rm
tr}_1(\overline{\xi}.\otimes \overline{\kappa})_k\mid_P\rightarrow
K(\overline{\eta},\overline{N})_k\otimes
p_P^*\overline{\kappa}\rightarrow 0,
\end{displaymath} we will get
\begin{displaymath}
({p_P}_g)_*[\widetilde{{\rm ch}}_g(\overline{\varepsilon}_k\otimes
p_P^*\overline{\kappa})]=({p_P}_g)_*[\widetilde{{\rm
ch}}_g(\overline{\varepsilon}_k)({p_P}_g)^*{\rm
ch}_g(\overline{\kappa})]=({p_p}_g)_*[\widetilde{{\rm
ch}}_g(\overline{\varepsilon}_k)]\cdot [{\rm
ch}_g(\overline{\kappa})].
\end{displaymath} Combing the two
computations above and the unique expression of $T$ via $C_T$, we
get the equality in the statement of this proposition.
\end{proof}

\begin{defn}\label{703}
An equivariant characteristic class $C$ for pairs of equivariant
hermitian vector bundles is said to be compatible with the
projection formula if it satisfies
\begin{displaymath}
C(\eta,N)=C(\mathcal{O}_Y,N)\cdot {\rm ch}_g(\eta).
\end{displaymath}
\end{defn}

The following is the main theorem in this section.

\begin{thm}\label{704}
A theory of equivariant singular Bott-Chern classes $T$ for
$\Sigma$ is compatible with the projection formula if and only if
the associated characteristic class $C_T$ is so.
\end{thm}
\begin{proof}
We first assume that $C_T$ is compatible with the projection
formula, then we compute
\begin{align*}
{i_g}_*C_T(\eta\otimes
i^*\kappa,N)&={i_g}_*(C_T(\mathcal{O}_Y,N)\cdot {\rm
ch}_g(\eta\otimes i^*\kappa))={i_g}_*(C_T(\mathcal{O}_Y,N)\cdot
{\rm ch}_g(\eta)\cdot{i_g}^*{\rm ch}_g(\kappa))\\
&={i_g}_*(C_T(\mathcal{O}_Y,N)\cdot {\rm ch}_g(\eta))\cdot{\rm
ch}_g(\kappa)={i_g}_*(C_T(\eta,N))\cdot{\rm ch}_g(\kappa).
\end{align*}
Therefore, by Proposition~\ref{702}, $T$ is compatible with the
projection formula.

For the other direction, assume that $T$ is compatible with the
projection formula. Using the definition of $C_T$, we compute
\begin{align*}
C_T(\eta,N)&=({\pi_P}_g)_*(T(K(\overline{\eta},\overline{N})))=({\pi_P}_g)_*(T(K(\overline{O}_Y,\overline{N})\otimes
\pi_p^*\overline{\eta}))\\
&=({\pi_P}_g)_*(T(K(\overline{O}_Y,\overline{N}))\cdot
{\pi_P}_g^*{\rm
ch}_g(\overline{\eta}))=({\pi_P}_g)_*(T(K(\overline{O}_Y,\overline{N})))\cdot{\rm
ch}_g(\overline{\eta})\\
&=C_T(\mathcal{O}_Y,N)\cdot{\rm ch}_g(\eta).
\end{align*} This
implies that $C_T$ is compatible with the projection formula.
\end{proof}

\section{Uniqueness of equivariant singular Bott-Chern classes}
Let $\Sigma$ be any set of equivariant hermitian embedded vector
bundles which satisfies the condition (Hui) and whose elements
have bounded equivariant ranks and bounded equivariant
codimensions. By Corollary~\ref{606}, it is possible to attach
$\Sigma$ a theory of equivariant singular Bott-Chern classes. But
unfortunately, such a theory is not unique. Our aim in this
section is to show that if we add another axiom to
Definition~\ref{505}, we will get a unique theory of equivariant
singular Bott-Chern classes for $\Sigma$ without the limitation of
the bounds of equivariant rank and codimension. Such a theory will
be called a theory of equivariant homogeneous singular Bott-Chern
classes. We shall also compare it with the theory of equivariant
singular Bott-Chern currents defined by J.-M. Bismut in \cite{Bi}.

Our startint point is again the Koszul construction. Let $Y$ be an
equivariant projective manifold. Assume that we are given two
equivariant hermitian vector bundles $\overline{\eta}$ and
$\overline{N}$ on $Y$. Let $P=\mathbb{P}(N\oplus \mathcal{O}_Y)$,
$P_0=\mathbb{P}(N_g\oplus\mathcal{O}_{Y_g})$ and let $i_\infty$ be
the zero section embedding. Suppose that $T$ is a theory of
equivariant singular Bott-Chern classes for $\Sigma$. Then by
definition, we have
\begin{displaymath}
{\rm dd}^cT(K(\overline{\eta},\overline{N}))=c_{{\rm
rk}{Q_g}}(\overline{Q}_g\mid_{P_0}){\rm
Td}_g^{-1}(\overline{Q}){\rm
ch}_g(\pi_P^*\overline{\eta})-(i_{\infty,g})_*({\rm
ch}_g(\overline{\eta}){\rm Td}_g^{-1}(\overline{N})).
\end{displaymath} Therefore, the class
\begin{displaymath}
\widetilde{e}_T(\overline{\eta},\overline{N}):=T(K(\overline{\eta},\overline{N}))\cdot{\rm
Td}_g(\overline{Q})\cdot{\rm ch}_g^{-1}({\pi_P}^*\overline{\eta})
\end{displaymath} satisfies the following differential equation
\begin{displaymath}
{\rm dd}^c\widetilde{e}_T(\overline{\eta},\overline{N})=c_{{\rm
rk}{Q_g}}(\overline{Q}_g\mid_{P_0})-\delta_{Y_g}.
\end{displaymath} Note
that by our descriptions in Proposition~\ref{405} (iii) and (iv),
the current $c_{{\rm
rk}{Q_g}}(\overline{Q}_g\mid_{P_0})-\delta_{Y_g}$ belongs to ${\rm
Im}(D^{{\rm rk}Q_g,{\rm rk}Q_g}(P_0)\hookrightarrow D^{{\rm
rk}Q_g,{\rm rk}Q_g}(P_g))$. Then it is natural to introduce the
following definition.

\begin{defn}\label{801}
Let $T$ be a theory of equivariant singular Bott-Chern classes for
some $\Sigma$. The class
$\widetilde{e}_T(\overline{\eta},\overline{N})$ is called the
Euler-Green class associated to $T$. We say that $T$ is
homogeneous if
\begin{displaymath}
\widetilde{e}_T(\overline{\eta},\overline{N})\in
\widetilde{\mathcal{U}}^{{\rm rk}Q_g-1,{\rm
rk}Q_g-1}(P_0):=D^{{\rm rk}Q_g-1,{\rm rk}Q_g-1}(P_0)/({{\rm
Im}\partial+{\rm Im}\overline{\partial}}) \end{displaymath} for
any element $(i: Y\rightarrow
X,\overline{N},\overline{\eta},\overline{\xi.})\in \Sigma$.
\end{defn}

\begin{rem}\label{802}
If $T$ is compatible with the projection formula, then the
Euler-Green class $\widetilde{e}_T(\overline{\eta},\overline{N})$
has nothing to do with the first variable.
\end{rem}

\begin{thm}\label{803}
Let $\Sigma$ be any set of equivariant hermitian embedded vector
bundles which satisfies the condition (Hui). Then there exists a
unique theory of equivariant homogeneous singular Bott-Chern
classes for $\Sigma$.
\end{thm}
\begin{proof}
We shall use a uniqueness theorem of the Euler-Green class in
non-equivariant case. That's the following.
\begin{lem}\label{804}
Let $i_\infty: Y\rightarrow P=\mathbb{P}(N\oplus\mathcal{O}_Y)$ be
a zero section embedding in non-equivariant setting. Denote
$D_\infty=\mathbb{P}(N)$. Then there exists a unique class
$\widetilde{e}(P,\overline{Q},i_\infty)\in
\widetilde{\mathcal{U}}^{{\rm rk}Q-1,{\rm rk}Q-1}(P)$ such that

(i). ${\rm dd}^c\widetilde{e}(P,\overline{Q},i_\infty)=c_{{\rm
rk}Q}(\overline{Q})-\delta_Y$;

(ii). $\widetilde{e}(P,\overline{Q},i_\infty)\mid_{D_\infty}=0$.
\end{lem}
We refer to \cite[Lemma 9.4]{BL} for a proof of this lemma. One
just need to pay attention to two points. Firstly, the restriction
isomorphism on analytic Deligne cohomology should be changed to
$H^{{\rm rk}Q-1,{\rm rk}Q-1}(P)\cong H^{{\rm rk}Q-1,{\rm
rk}Q-1}(D_\infty)$ on Dolbeault cohomology which can be deduced
from the classical projective bundle theorem for deRham cohomology
and Hodge decomposition. Secondly, the existence of a preimage of
$c_{{\rm rk}Q}(\overline{Q})-\delta_Y$ under ${\rm dd}^c$ is a
consequence of \cite[Theorem 1.2.1]{GS2}. What we want to indicate
is that this lemma naturally leads to a similar result in the
equivariant setting by using Proposition~\ref{405} (iii) and (iv).
The result reads: there exists a unique class
$\widetilde{e}(P,\overline{Q},i_\infty)\in {\rm
Im}(\widetilde{\mathcal{U}}^{{\rm rk}Q_g-1,{\rm
rk}Q_g-1}(P_0)\hookrightarrow \widetilde{\mathcal{U}}^{{\rm
rk}Q_g-1,{\rm rk}Q_g-1}(P_g))$ such that ${\rm
dd}^c\widetilde{e}(P,\overline{Q},i_\infty)=c_{{\rm
rk}Q_g}(\overline{Q}_g\mid_{P_0})-\delta_{Y_g}$ and
$\widetilde{e}(P,\overline{Q},i_\infty)\mid_{D_{\infty,g}}=0$.
Moreover, by convention, we shall identify $D^{-1,-1}(P_g)$ with
the zero space.

Now assume that $T$ is a theory of equivariant homogeneous
singular Bott-Chern classes. Since the restriction of the Koszul
resolution $K(\overline{\eta},\overline{N})$ to $D_{\infty,g}$ is
equivariantly and orthogonally split, then we have
$T(K(\overline{\eta},\overline{N}))\mid_{D_{\infty,g}}=0$. Thus
the restriction of the Euler-Green class
$\widetilde{e}_T(\overline{\eta},\overline{N})$ to $D_{\infty,g}$
is equal to $0$ by definition. Therefore, by uniqueness, we get
$\widetilde{e}_T(\overline{\eta},\overline{N})=\widetilde{e}(P,\overline{Q},i_\infty)$
and hence
\begin{displaymath}
T(K(\overline{\eta},\overline{N}))=\widetilde{e}(P,\overline{Q},i_\infty)\cdot
{\rm Td}_g^{-1}(\overline{Q})\cdot {\rm
ch}_g(\pi_P^*\overline{\eta}).
\end{displaymath} This equation
implies that the characteristic class
\begin{displaymath}
C_T(\eta,N)=({\pi_P}_g)_*T(K(\overline{\eta},\overline{N}))
\end{displaymath} is independent of the theory $T$. So the uniqueness of $T$
follows from Theorem~\ref{601}.

For the existence, we define
\begin{displaymath}
C(\overline{\eta},\overline{N})=({\pi_P}_g)_*(\widetilde{e}(P,\overline{Q},i_\infty)\cdot{\rm
Td}_g^{-1}(\overline{Q})\cdot{\rm ch}_g(\pi_P^*\overline{\eta})).
\end{displaymath} By the differential equation that
$\widetilde{e}(P,\overline{Q},i_\infty)$ satisfies, one can easily
prove that $C(\overline{\eta},\overline{N})$ is ${\rm
dd}^c$-closed. This is the only important point for us to use the
same principle as in the proof of Theorem~\ref{602} to define a
theory of equivariant homogeneous singular Bott-Chern classes $T$
such that $C_T=C$. The last equality says that
$C(\overline{\eta},\overline{N})$ is actually independent of the
choices of the metrics since $C_T$ is so. Thus we can just write
\begin{displaymath}
C(\eta,N)=({\pi_P}_g)_*(\widetilde{e}(P,\overline{Q},i_\infty)\cdot{\rm
Td}_g^{-1}(\overline{Q})\cdot{\rm ch}_g(\pi_P^*\overline{\eta})).
\end{displaymath} And this is compatible with Theorem~\ref{602}.
\end{proof}

Since the class $\widetilde{e}(P,\overline{Q},i_\infty)$ has
nothing to do with the vector bundle $\eta$, the following remark
looks more natural.

\begin{rem}\label{805}
If $T$ is compatible with the projection formula, then $T$ is
homogeneous if and only if
$\widetilde{e}_T(\overline{\mathcal{O}_Y},\overline{N})=\widetilde{e}(P,\overline{Q},i_\infty)$.
\end{rem}

We reformulate Theorem~\ref{803} in an axiomatical way.

\begin{thm}\label{806}
There exists a unique way to associate to each equivariant
hermitian embedded vector bundle $\overline{\Xi}=(i: Y\rightarrow
X,\overline{N},\overline{\eta},\overline{\xi}.)$ a class of
currents
\begin{displaymath}
T^h(\overline{\Xi})\in \widetilde{\mathcal{U}}(X_g,N_{g,0}^\vee)
\end{displaymath}
which we call equivariant homogeneous singular Bott-Chern class,
satisfying the following properties

(i). (Differential equation) The following equality holds
\begin{displaymath}
{\rm dd}^cT^h(\overline{\Xi})=\sum_{j}(-1)^j[{\rm
ch}_g(\overline{\xi}_j)]-{i_g}_*([{\rm ch}_g(\overline{\eta}){\rm
Td}_g^{-1}(\overline{N})]).
\end{displaymath}

(ii). (Functoriality) For every equivariant morphism $f:
X'\rightarrow X$ of projective manifolds which is transversal to
$Y$, we have
\begin{displaymath}
f_g^*T^h(\overline{\Xi})=T^h(f^*\overline{\Xi}).
\end{displaymath}

(iii). (Normalization) Let $\overline{A}.$ be an equivariantly and
orthogonally split exact sequence of equivariant hermitian vector
bundles. Write $\overline{\Xi}\oplus
\overline{A}.=(i,\overline{N},\overline{\eta},\overline{\xi}.\oplus\overline{A}.)$.
Then $T^h(\overline{\Xi})=T^h(\overline{\Xi}\oplus\overline{A}.)$.
Moreover, if $X={\rm Spec}(\C)$ is one point, $Y=\emptyset$ and
$\overline{\xi}.=0$, then $T^h(\overline{\Xi})=0$.

(iv). (Homogeneity) For any Koszul construction, we have
\begin{displaymath}
T^h(K(\overline{\eta},\overline{N}))\cdot {\rm
Td}_g(\overline{Q})\cdot {\rm ch}_g(\pi_P^*\overline{\eta})\in
{\rm Im}(\widetilde{\mathcal{U}}^{{\rm rk}Q_g-1,{\rm
rk}Q_g-1}(P_0)\hookrightarrow \widetilde{\mathcal{U}}^{{\rm
rk}Q_g-1,{\rm rk}Q_g-1}(P_g)).
\end{displaymath}
\end{thm}

\begin{prop}\label{807}
The theory of equivariant homogeneous singular Bott-Chern classes
is compatible with the projection formula.
\end{prop}
\begin{proof}
By definition, we compute
\begin{align*}
C_{T^h}(\eta,N)&=({\pi_P}_g)_*T^h(K(\overline{\eta},\overline{N}))=({\pi_P}_g)_*(\widetilde{e}(P,\overline{Q},i_\infty)\cdot{\rm
Td}_g^{-1}(\overline{Q})\cdot{\rm
ch}_g(\pi_P^*\overline{\eta}))\\
&=({\pi_P}_g)_*(\widetilde{e}(P,\overline{Q},i_\infty)\cdot{\rm
Td}_g^{-1}(\overline{Q}))\cdot{\rm
ch}_g(\overline{\eta})=C_{T^h}(\mathcal{O}_Y,N)\cdot{\rm
ch}_g(\eta). \end{align*} Then this proposition follows from
Theorem~\ref{704}
\end{proof}

The equivariant and non-equivariant homogeneous singular
Bott-Chern classes are related by the following proposition.

\begin{prop}\label{808}
Let $S^h$ be the non-equivariant homogeneous singular Bott-Chern
classes defined in \cite{BL}. Assume that $\overline{\eta}$ is an
equivariant hermitian vector bundle whose restriction to the fixed
point submanifold has no non-zero degree part. Then we have
\begin{displaymath}
T^h(K(\overline{\eta},\overline{N}))\cdot {\rm
Td}_g(\overline{Q})=S^h(K(\overline{\eta}_g,\overline{N}_g))\cdot{\rm
Td}(\overline{Q}_g).
\end{displaymath}
\end{prop}
\begin{proof}
We first suppose that $\overline{\eta}$ is the trivial bundle
$\overline{\mathcal{O}_Y}$ equipped with the trivial
$g$-structure. By equation \cite[(9.8)]{BL}, we have
\begin{displaymath}
S^h(K(\overline{\mathcal{O}_{Y_g}},\overline{N}_g))=\widetilde{e}(P_0,\overline{Q}_g,i_{\infty,0})\cdot{\rm
Td}^{-1}(\overline{Q}_g) \end{displaymath} where $i_{\infty,0}$ is
the zero section embedding from $Y_g$ to
$P_0=\mathbb{P}(N_g\oplus\mathcal{O}_{Y_g})$. Note that by the
definition of homogeneity in our paper, the class
$\widetilde{e}(P,\overline{Q},i_\infty)$ is equal to
$\widetilde{e}(P_0,\overline{Q}_g,i_{\infty,0})$ in ${\rm
Im}(\widetilde{\mathcal{U}}^{{\rm rk}Q_g-1,{\rm
rk}Q_g-1}(P_0)\hookrightarrow \widetilde{\mathcal{U}}^{{\rm
rk}Q_g-1,{\rm rk}Q_g-1}(P_g))$. This implies that
\begin{displaymath}
T^h(K(\overline{\mathcal{O}_Y},\overline{N}))\cdot {\rm
Td}_g(\overline{Q})=S^h(K(\overline{\mathcal{O}_{Y_g}},\overline{N}_g))\cdot{\rm
Td}(\overline{Q}_g). \end{displaymath} In general case, since the
restriction of $\overline{\eta}$ to $Y_g$ is supposed to have no
non-zero degree part, we have ${\rm
ch}_g(\pi_P^*\overline{\eta})={\pi_P}_g^*{\rm
ch}(\overline{\eta}_g)$. Moreover, the class
$S^h(K(\overline{\mathcal{O}_{Y_g}},\overline{N}_g))\cdot{\rm
Td}(\overline{Q}_g)$ belongs to $D^{{\rm rk}Q_g-1,{\rm
rk}Q_g-1}(P_0)/({\rm Im}\partial+{\rm Im}\overline{\partial})$, we
then can compute
\begin{displaymath}
T^h(K(\overline{\mathcal{O}_Y},\overline{N}))\cdot {\rm
Td}_g(\overline{Q})\cdot {\rm
ch}_g(\pi_P^*\overline{\eta})=S^h(K(\overline{\mathcal{O}_{Y_g}},\overline{N}_g))\cdot{\rm
Td}(\overline{Q}_g)\cdot \pi_{P_0}^*{\rm ch}(\overline{\eta}_g).
\end{displaymath}
This equality implies that
$T^h(K(\overline{\eta},\overline{N}))\cdot {\rm
Td}_g(\overline{Q})=S^h(K(\overline{\eta}_g,\overline{N}_g))\cdot{\rm
Td}(\overline{Q}_g)$ because $T^h$ and $S^h$ are both compatible
with the projection formula.
\end{proof}

In general, let $X$ be a complex manifold and let $\overline{E}$
be a hermitian holomorphic vector bundle of rank $r$ on $X$.
Assume that $s$ is a holomorphic section of $E$ which is
transversal to the zero section. Denote by $Y$ the zero locus of
$s$. In \cite[Proposition 9.13]{BL}, the authors have shown that
there is a unique way to attach to each $(X,\overline{E},s)$ a
class of currents $\widetilde{e}(X,\overline{E},s)\in
\widetilde{\mathcal{U}}^{r-1,r-1}(X,N^\vee_{Y,0})$ which satisfies
some axiomatic properties. Such class was also constructed by
J.-M. Bismut, H. Gillet and C. Soul\'{e} in \cite{BGS2}. We shall
use this fact to generalize Proposition~\ref{808} in the following
way. Assume that all notations above are $g$-equivariant, then
there is a global equivariant Koszul resolution
\begin{displaymath}
K(\overline{E}):\quad 0\rightarrow
\wedge^r\overline{E}^\vee\rightarrow\cdots\rightarrow
\overline{E}^\vee\rightarrow\overline{\mathcal{O}_X}\rightarrow
i_*\overline{\mathcal{O}_Y}\rightarrow 0. \end{displaymath} So we
get an equivariant hermitian embedded vector bundle
$(i,\overline{N}_{X/Y},\overline{\mathcal{O}_Y},K(\overline{E}))$
such that $\overline{N}_{X/Y}$ is isometric to $i^*\overline{E}$.
One can carry out the proof of \cite[Prop. 9.18]{BL} word by word
(adding subscript $g$) to prove the following equality
\begin{displaymath}
T^h(i,\overline{N}_{X/Y},\overline{\mathcal{O}_Y},K(\overline{E}))=\widetilde{e}(X_g,\overline{E}_g,s_g)\cdot{\rm
Td}_g^{-1}(\overline{E}). \end{displaymath} This equality and
\cite[Prop. 9.18]{BL} imply the following result.

\begin{prop}\label{809}
Let notations and assumptions be as above, then we have
\begin{displaymath}
T^h(i,\overline{N}_{X/Y},\overline{\mathcal{O}_Y},K(\overline{E}))\cdot{\rm
Td}_g(\overline{E})=S^h(i_g,\overline{N}_{{X_g}/{Y_g}},\overline{\mathcal{O}_{Y_g}},K(\overline{E}_g))\cdot{\rm
Td}(\overline{E}_g).
\end{displaymath}
\end{prop}

We now recall the construction of the equivariant Bott-Chern
singular currents given by J.-M. Bismut in \cite{Bi}. This
construction was realized via some current valued zeta function
which involves the supertraces of Quillen's superconnections. We
would like to indicate that Bismut's singular current defines a
class which agrees with our definition of equivariant singular
Bott-Chern class only in some certain situation. Nevertheless, it
is easy to use Bismut's results to define a theory of equivariant
singular Bott-Chern classes in the sense of Definition~\ref{505}.
We shall prove that such a theory is homogeneous.

Let $i: Y\rightarrow X$ be a closed immersion of equivariant
projective manifolds, and let
$\overline{\Xi}=(i,\overline{N},\overline{\eta},\overline{\xi}.)$
be an equivariant hermitian embedded vector bundle. We denote the
differential of the complex $\xi.$ by $v$. Note that $\xi.$ is
acyclic outside $Y$ and the homology sheaves of its restriction to
$Y$ are locally free. We write $H_n=\mathcal{H}_n(\xi.\mid_Y)$ and
define a $\Z$-graded bundle $H=\bigoplus_nH_n$. For $y\in Y$ and
$u\in TX_y$, we denote by $\partial_uv(y)$ the derivative of $v$
at $y$ in the direction $u$ in any given holomorphic
trivialization of $\xi.$ near $y$. Then the map $\partial_uv(y)$
acts on $H_y$ as a chain map, and this action only depends on the
image $z$ of $u$ in $N_y$. So we get a chain complex of
holomorphic vector bundles $(H,\partial_zv)$.

Let $\pi$ be the projection from the normal bundle $N$ to $Y$,
then we have a canonical identification of $\Z$-graded chain
complexes
\begin{displaymath}
(\pi^*H,\partial_zv)\cong(\pi^*(\wedge^\bullet
N^\vee\otimes\eta),\sqrt{-1}i_z).
\end{displaymath} Moreover, such
an identification is an identification of $g$-bundles. By finite
dimensional Hodge theory, for each $y\in Y$, there is a canonical
isomorphism
\begin{displaymath}
H_y\cong\{f\in \xi._y\mid vf=0, v^*f=0\} \end{displaymath} where
$v^*$ is the dual of $v$ with respect to the metrics on $\xi.$.
This means that $H$ can be regarded as a smooth $\Z$-graded
$g$-equivariant subbundle of $\xi$ so that it carries an induced
$g$-invariant metric. On the other hand, we endow $\wedge^\bullet
N^\vee\otimes \eta$ with the metric induced from $\overline{N}$
and $\overline{\eta}$.

\begin{defn}\label{810}
We say that the metrics on the complex of equivariant hermitian
vector bundles $\overline{\xi}.$ satisfy Bismut assumption (A) if
the identification $(\pi^*H,\partial_zv)\cong(\pi^*(\wedge^\bullet
N^\vee\otimes\eta),\sqrt{-1}i_z)$ also identifies the metrics.
\end{defn}

\begin{prop}\label{811}
There always exist $g$-invariant metrics on $\xi.$ which satisfy
Bismut assumption (A) with respect to $\overline{N}$ and
$\overline{\eta}$.
\end{prop}
\begin{proof}
This is \cite[Proposition 3.5]{Bi}.
\end{proof}

Let $\nabla^{\xi}$ be the canonical hermitian holomorphic
connection on $\xi.$, then for $u>0$, we may define a
$g$-invariant superconnection
\begin{displaymath}
C_u:=\nabla^\xi+\sqrt{u}(v+v^*)
\end{displaymath} on the
$\Z_2$-graded vector bundle $\xi$. Let $\Phi$ be the map
$\alpha\in \wedge(T_\R^*X_g)\rightarrow (2\pi i)^{-{\rm
deg}\alpha/2}\alpha\in \wedge(T_\R^*X_g)$ and denote
\begin{displaymath}
({\rm Td}_g^{-1})'(\overline{N}):=\frac{\partial}{\partial
b}\mid_{b=0}({\rm Td}_g(b\cdot {\rm
Id}-\frac{\Omega^{\overline{N}}}{2\pi i})^{-1}).
\end{displaymath}

\begin{lem}\label{812}
Let $N_H$ be the number operator on the complex $\xi.$ i.e. it
acts on $\xi_j$ as multiplication by $j$, then for $s\in \C$ and
$0< {\rm Re}(s)<\frac{1}{2}$, the current valued zeta function
\begin{displaymath}
Z_g(\overline{\xi}.)(s):=\frac{1}{\Gamma(s)}\int_0^\infty
u^{s-1}[\Phi{\rm Tr_s}(N_Hg{\rm exp}(-C_u^2))+({\rm
Td}_g^{-1})'(\overline{N}){\rm
ch}_g(\overline{\eta})\delta_{Y_g}]{\rm d}u \end{displaymath} is
well-defined on $X_g$ and it has a meromorphic continuation to the
complex plane which is holomorphic at $s=0$.
\end{lem}

\begin{defn}\label{813}
The equivariant Bott-Chern singular current on $X_g$ associated to
the resolution $\overline{\xi}.$ is defined as
\begin{displaymath}
T_g(\overline{\xi}.):=\frac{\partial}{\partial
s}\mid_{s=0}Z_g(\overline{\xi}.)(s).
\end{displaymath}
\end{defn}

\begin{thm}\label{814}
The current $T_g(\overline{\xi}.)$ is a sum of $(p,p)$-currents
and it satisfies the differential equation
\begin{displaymath}
{\rm dd}^cT_g(\overline{\xi}.)={i_g}_*{\rm
ch}_g(\overline{\eta}){\rm
Td}_g^{-1}(\overline{N})-\sum_k(-1)^k{\rm ch}_g(\overline{\xi}_k).
\end{displaymath}
Moreover, the wave front set of $T_g(\overline{\xi}.)$ is
contained in $N^\vee_{g,0}$.
\end{thm}

For any equivariant hermitian embedded vector bundle
$\overline{\Xi}_0=(i,\overline{N},\overline{\eta},(\xi.,h_0^\xi))$,
we may construct a new embedded bundle
$\overline{\Xi}_1=(i,\overline{N},\overline{\eta},(\xi.,h_1^\xi))$
such that the metrics $h_1^\xi$ satisfies Bismut assumption (A).
Then we may attach to $\overline{\Xi}_0$ an element in
$\widetilde{\mathcal{U}}(X_g)$ defined as
\begin{displaymath}
T^B(\overline{\Xi}_0)=-T_g(\xi.,h_1^\xi)+\sum_k(-1)^k\widetilde{{\rm
ch}}_g(\xi_k,h_0^{\xi_k},h_1^{\xi_k}).
\end{displaymath}

\begin{thm}\label{815}
The assignment that, to each equivariant hermitian embedded vector
bundle $\overline{\Xi}_0$, associates the current
$T^B(\overline{\Xi}_0)$, is a theory of equivariant homogeneous
singular Bott-Chern classes.
\end{thm}
\begin{proof}
We first show that $T^B(\overline{\Xi}_0)$ is well-defined.
Actually, let
$\overline{\Xi}_2=(i,\overline{N},\overline{\eta},(\xi.,h_2^\xi))$
be another embedded bundle such that the metrics $h_2^\xi$ satisfy
Bismut assumption (A), then by \cite[Theorem 3.14]{KR1} we have
\begin{displaymath}
T_g(\overline{\xi}_1)-T_g(\overline{\xi}_2)=-\sum_k(-1)^k\widetilde{{\rm
ch}}_g(\xi_k,h_1^{\xi_k},h_2^{\xi_k}). \end{displaymath} Note that
we have the equality
\begin{displaymath}
\widetilde{{\rm
ch}}_g(\xi_k,h_0^{\xi_k},h_1^{\xi_k})+\widetilde{{\rm
ch}}_g(\xi_k,h_1^{\xi_k},h_2^{\xi_k})+\widetilde{{\rm
ch}}_g(\xi_k,h_2^{\xi_k},h_0^{\xi_k})=0. \end{displaymath} So we
obtain that $T^B(\overline{\Xi}_0)$ does not depend on the choice
of the metrics which satisfy Bismut assumption (A) and hence it is
well-defined.

Secondly, the fact that the equivariant singular current
$T^B(\overline{\Xi}_0)$ satisfies the differential equation in
Definition~\ref{505} follows from Theorem~\ref{814} and the
definition of $\widetilde{{\rm ch}}_g$.

The functoriality property for $T^B(\overline{\Xi}_0)$ follows
from the same property for $T_g$ and for $\widetilde{{\rm ch}}_g$.

For the normalization property, let $\overline{A}.$ be an
equivariantly and orthogonally split exact sequence of equivariant
hermitian vector bundles, then using \cite[Theorem 3.14]{KR1}
again we have
\begin{displaymath}
T_g(\overline{\xi}.\oplus\overline{A}.)=T_g(\overline{\xi}.)+T_g(\overline{A}.).
\end{displaymath}
By \cite[Corollary 3.10]{KR1}, if $\overline{A}.$ is equivariantly
and orthogonally split, then $T_g(\overline{A}.)$ is equal to
zero. So by definition we finally get
$T^B(\overline{\Xi}\oplus\overline{A}.)=T^B(\overline{\Xi})$.

At last, by \cite[Lemma 3.15]{KR1}, with the hypothesis before
Proposition~\ref{809} we have the following equality
\begin{displaymath}
T^B(i,\overline{N}_{X/Y},\overline{\mathcal{O}_Y},K(\overline{E}))=\widetilde{e}(X_g,\overline{E}_g,s_g)\cdot
{\rm
Td}_g^{-1}(\overline{E})=T^h(i,\overline{N}_{X/Y},\overline{\mathcal{O}_Y},K(\overline{E})).
\end{displaymath}
Since $T^B$ and $T^h$ are both compatible with the projection
formula, this equality implies that $C_{T^B}=C_{T^h}$ and hence
$T^B=T^h$ by Theorem~\ref{601}. So $T^B$ is homogeneous which
completes the whole proof.
\end{proof}

\section{Concentration formula}
In the last section, we shall prove a concentration formula for
equivariant homogeneous singular Bott-Chern class. We call it
concentration formula because it can be used to prove a statement
which generalizes the concentration theorem in algebraic
$K$-theory (cf. \cite{Th}) to the context of Arakelov geometry. We
deal with this in another paper. Before describing the
concentration formula, we introduce some basic concepts.

\begin{defn}\label{901}
Let $X$ be a complex manifold and let $\overline{\xi}.$ be a
bounded complex of hermitian vector bundles on $X$. We say
$\overline{\xi}.$ is standard if the homology sheaves of
$\overline{\xi}.$ are all locally free and they are endowed with
some hermitian metrics. We shall write a standard complex as
$(\overline{\xi}.,h^H)$ to emphasize the choice of the metrics on
homology sheaves.
\end{defn}

Now let $X$ be a $\mu_n$-equivariant projective manifold, we
consider a special closed immersion $i: X_g\hookrightarrow X$. For
an equivariant hermitian embedded vector bundle
$\overline{\Xi}=(i,\overline{N},\overline{\eta},\overline{\xi}.)$,
we always assume that the metrics on $\xi.$ satisfy Bismut
assumption (A). In this case, the restriction of $\overline{\xi}.$
to $X_g$ is a standard complex according to our discussion in last
section, the metrics on homology bundles are induced by the
metrics on $\xi.\mid_{X_g}$. Note that we can split
$\overline{\xi}.\mid_{X_g}$ into a series of short exact sequences
\begin{displaymath}
0\rightarrow\overline{{\rm Im}}\rightarrow\overline{{\rm
Ker}}\rightarrow\wedge^\bullet\overline{N}^\vee\otimes\overline{\eta}\rightarrow
0 \end{displaymath} and
\begin{displaymath}
0\rightarrow\overline{{\rm
Ker}}\rightarrow\overline{\xi}.\mid_{X_g}\rightarrow\overline{{\rm
Im}}\rightarrow 0. \end{displaymath} Denote the alternating sum of
the equivariant secondary characteristic classes of the short
exact sequences above by $\widetilde{{\rm
ch}}_g(\overline{\xi}.,h^H)$ such that it satisfies the following
differential equation
\begin{displaymath}
{\rm dd}^c\widetilde{{\rm ch}}_g(\overline{\xi}.,h^H)={\rm
ch}_g(\overline{\eta}){\rm
Td}_g^{-1}(\overline{N})-\sum_j(-1)^j{\rm ch}_g(\overline{\xi}_j).
\end{displaymath}
With this observation, we can introduce the following proposition.

\begin{prop}\label{additional}
Let $\overline{\chi}:\quad 0\rightarrow
\overline{\eta}_n\rightarrow\cdots\rightarrow\overline{\eta}_1\rightarrow\overline{\eta}_0\rightarrow0$
be an exact sequence of equivariant hermitian vector bundles on
$X_g$, and let $\overline{\varepsilon}: 0\rightarrow
\overline{\xi}_{n,\cdot}\rightarrow
\cdots\rightarrow\overline{\xi}_{1,\cdot}\rightarrow\overline{\xi}_{0,\cdot}\rightarrow
0$ be an exact sequence of resolutions of $i_*\overline{\chi}$ on
$X$. As usual we write $\overline{\varepsilon}_k$ for the exact
sequence
\begin{displaymath}
0\rightarrow \overline{\xi}_{n,k}\rightarrow\cdots\rightarrow
\overline{\xi}_{1,k}\rightarrow\overline{\xi}_{0,k}\rightarrow 0.
\end{displaymath}
Then we have the following equality in $\widetilde{A}(X_g)$
\begin{displaymath}
\sum_{j=0}^n(-1)^j\widetilde{{\rm
ch}}_g(\overline{\xi}_{j,\cdot},h^H)=\widetilde{{\rm
ch}}_g(\overline{\chi}){\rm
Td}_g^{-1}(\overline{N})-\sum_k(-1)^k\widetilde{{\rm
ch}}_g(\overline{\varepsilon}_k).
\end{displaymath}
\end{prop}
\begin{proof}
Note that the fixed point submanifold of $X\times \mathbb{P}^1$ is
exactly $X_g\times\mathbb{P}^1$, we know that the construction of
the first transgression exact sequence is compatible with
restriction to the fixed point submanifold. This means ${\rm
tr}_1(\overline{\varepsilon}.)_j\mid_{X_g\times\mathbb{P}^1}$ is
equal to ${\rm tr}_1(\overline{\varepsilon}.\mid_{X_g})_j$.
Therefore, one can use the same approach as in the proof of
Proposition~\ref{507} to verify the equality in the statement of
this proposition.
\end{proof}

\begin{thm}\label{902}(Concentration formula)
Let notations and assumptions be as above. Assume that
$\overline{\Xi}=(i: X_g\rightarrow
X,\overline{N},\overline{\eta},\overline{\xi}.)$ is an equivariant
hermitian embedded vector bundle such that the metrics on $\xi.$
satisfy Bismut assumption (A). Then in $\widetilde{A}(X_g)$, we
have the equality
\begin{displaymath}
T^h(\overline{\Xi})=-\widetilde{{\rm ch}}_g(\overline{\xi}.,h^H).
\end{displaymath}
\end{thm}

Before proving this theorem, we first investigate the problem for
a simple case where the hypothesis is the same as before
Proposition~\ref{809}. That means there exists an equivariant
hermitian vector bundle $\overline{E}$ on $X$ which admits a
$g-$invariant regular section $s$ such that $X_g$ is the zero
locus of $s$ and $i^*\overline{E}$ is isometric to
$\overline{N}_{X/{X_g}}$. From this we know that $\overline{E}_g$
is the zero bundle, so
$S^h(i_g,\overline{N}_{X_g/{X_g}},\overline{\mathcal{O}_{X_g}},K(\overline{E}_g))=0$
and hence
$T^h(i,\overline{N}_{X/{X_g}},\overline{\mathcal{O}_{X_g}},K(\overline{E}))=0$
by Proposition~\ref{809}. On the other hand, $\widetilde{{\rm
ch}}_g(\wedge^\bullet\overline{E},h^H)$ is definitely equal to $0$
since the metrics $h^H$ are supposed to be induced from
$\overline{E}$ and these metrics satisfy Bismut assumption (A).
Therefore the concentration formula is trivially true for this
case.

\begin{proof}(of Theorem~\ref{902})
We use the same notations as in the proof of Theorem~\ref{601},
then we have the following expression
\begin{align*}
T^h(\overline{\Xi})&=-({p_W}_g)_*(U\cdot(\sum_k(-1)^k{\rm
ch}_g({\rm tr}_1(\overline{\xi}.)_k)-{j_g}_*({\rm
ch}_g(p_{X_g}^*\overline{\eta}){\rm
Td}_g^{-1}(\overline{N}'))))\\
&\qquad\qquad\qquad\qquad\qquad-\sum_k(-1)^k({p_P}_g)_*[\widetilde{{\rm
ch}}_g(\overline{\varepsilon}_k)]+C_{T^h}(\eta,N). \end{align*}
Since $i_g$ is the identity map, we know that the deformation to
the normal cone $W(i_g)$ is equal to $X_g\times \mathbb{P}^1$.
Moreover $W(i_g)$ is a disjoint union of some connected components
of $W_g$ and the map $j_g$ factors through $W(i_g)$. We shall
write $W_0$ for $W(i_g)$ for simplicity and we shall denote by
$W_\bot$ the other components of $W_g$. Now we restrict the sum
\begin{displaymath}
L:=\sum_k(-1)^k{\rm ch}_g({\rm
tr}_1(\overline{\xi}.)_k)-{j_g}_*({\rm
ch}_g(p_{X_g}^*\overline{\eta}){\rm Td}_g^{-1}(\overline{N}'))
\end{displaymath} to $W_\bot$ and $W_0$. Over $W_\bot$ we get
$L\mid_{W_\bot}=\sum_k(-1)^k{\rm ch}_g({\rm
tr}_1(\overline{\xi}.)_k\mid_{W_\bot})$ which can be rewritten as
${\rm dd}^c\widetilde{{\rm ch}}_g({\rm
tr}_1(\overline{\xi}.)\mid_{W_\bot})$. Similarly, over $W_0$ we
get
\begin{displaymath}
L\mid_{W_0}=\sum_k(-1)^k{\rm ch}_g({\rm
tr}_1(\overline{\xi}.)_k\mid_{W_0})-({\rm
ch}_g(p_{X_g}^*\overline{\eta}){\rm Td}_g^{-1}(\overline{N}'))
\end{displaymath}
which can be rewritten as $-{\rm dd}^c\widetilde{{\rm ch}}_g({\rm
tr}_1(\overline{\xi}.)\mid_{W_0},h^H)$ since in this case ${\rm
tr}_1(\overline{\xi}.)\mid_{W_0}$ is clearly a standard complex in
the sense of Definition~\ref{901}. Moreover, a totally similar
argument to the observation given before this proof shows that
$T^h(K(\overline{\eta},\overline{N}))$ is equal to $0$ so that
$C_{T^h}(\eta,N)$ is equal to $0$. Furthermore, the exact sequence
$K(\overline{\eta},\overline{N})\mid_{W_\bot\cap P}$ is
equivariantly and orthogonally split. Therefore, by
Remark~\ref{209}, we get
\begin{displaymath}
-\widetilde{{\rm ch}}_g({\rm
tr}_1(\overline{\xi}.)\mid_{W_\bot\cap
P})=\sum_k(-1)^k\widetilde{{\rm
ch}}_g(\overline{\varepsilon}_k)\mid_{W_\bot}. \end{displaymath}
This means
\begin{align*}
U\cdot L\mid_{W_\bot}&=U\cdot {\rm dd}^c\widetilde{{\rm
ch}}_g({\rm tr}_1(\overline{\xi}.)\mid_{W_\bot})={\rm
dd}^cU\cdot\widetilde{{\rm ch}}_g({\rm
tr}_1(\overline{\xi}.)\mid_{W_\bot})\\
&=\widetilde{{\rm ch}}_g({\rm
tr}_1(\overline{\xi}.)\mid_{W_\bot\cap
P})=-\sum_k(-1)^k\widetilde{{\rm
ch}}_g(\overline{\varepsilon}_k)\mid_{W_\bot}. \end{align*}
Combing these computations above, we may reformulate
$T^h(\overline{\Xi})$ as
\begin{displaymath}
T^h(\overline{\Xi})=-(p_{W_0})_*(U\cdot
L\mid_{W_0})-\sum_k(-1)^k(p_{P_0})_*[\widetilde{{\rm
ch}}_g(\overline{\varepsilon}_k)\mid_{P_0}].
\end{displaymath}
Similar to ${\rm tr}_1(\overline{\xi}.)\mid_{W_0}$,
$K(\overline{\eta},\overline{N})\mid_{P_0}$ is also a standard
complex. Since the metrics on the Koszul resolution are supposed
to satisfy Bismut assumption (A), we know that $\widetilde{{\rm
ch}}_g(K(\overline{\eta},\overline{N})\mid_{P_0},h^H)$ is equal to
$0$. Then by Proposition~\ref{additional}, we have that
\begin{displaymath}
\widetilde{{\rm ch}}_g({\rm tr}_1(\overline{\xi}.)\mid_{W_0\cap
P},h^H)=\sum_k(-1)^k\widetilde{{\rm
ch}}_g(\overline{\varepsilon}_k\mid_{P_0})=\sum_k(-1)^k\widetilde{{\rm
ch}}_g(\overline{\varepsilon}_k)\mid_{P_0}. \end{displaymath}
Together with the fact that ${\rm
tr}_1(\overline{\xi}.)\mid_{X_g\times \{0\}}$ is isometric to
$\overline{\xi}.\mid_{X_g}$, we finally get
\begin{displaymath}
T^h(\overline{\Xi})=-\widetilde{{\rm ch}}_g({\rm
tr}_1(\overline{\xi}.)\mid_{W_0},h^H)\mid_{X_g\times\{0\}}=-\widetilde{{\rm
ch}}_g(\overline{\xi}.,h^H)
\end{displaymath} which completes the proof.
\end{proof}

\hspace{5cm} \hrulefill\hspace{5.5cm}

D\'{e}partement de Math\'{e}matiques, B\^{a}timent 425,
Universit\'{e} Paris-Sud, 91405 Orsay cedex, France

E-mail: shun.tang@math.u-psud.fr


\begin{thebibliography}{200}
\bibitem{BGS1}
J.-M. Bismut, H. Gillet and C. Soul\'{e}, \emph{Analytic torsion
and holomorphic determinant bundles I}, Comm. Math. Phys.
\textbf{115}(1988), 49-78.

\bibitem{BGS2}
J.-M. Bismut, H. Gillet and C. Soul\'{e}, Complex immersions and
Arakelov geometry, in \emph{Grothendieck Festschrift I}, P.
Cartier and al.(eds.), Birkha\"{u}ser, 1990.

\bibitem{BGS3}
J.-M. Bismut, H. Gillet and C. Soul\'{e}, \emph{Bott-Chern current
and complex immersions}, Duke Math. J. \textbf{60}(1990), 255-284.

\bibitem{Bi}
J.-M. Bismut, \emph{Equivariant immersions and Quillen metrics},
J. Differential Geom. \textbf{41}(1995), 53-157.

\bibitem{BL}
J. I. Burgos Gil, R. Li\c{t}canu, \emph{Singular Bott-Chern
classes and the arithmetic Grothendieck Riemann Roch theroem for
closed immersions}, Documenta Math. \textbf{15}(2010), 73-176,
available online: http://www.math.uiuc.edu/documenta/

\bibitem{Ei}
S. Eilenberg, \emph{Homological dimension and local syzygies},
Annals of Math. \textbf{64}(1956), 328-336.

\bibitem{GBI}
A. Grothendieck, P. Berthelot and L. Illusie, SGA6,
\emph{Th\'{e}orie des intersections et th\'{e}or\`{e}me de
Riemann-Roch}, Lecture Notes in Maththematics \textbf{225},
Springer-Verlag, Berlin-Heidelberg-New York, 1971.

\bibitem{GH}
P. Griffiths and J. Harris, \emph{Principles of Algebraic
Geometry}, John Wiley and Sons, 1978.

\bibitem{GS2}
H. Gillet and C. Soul\'{e}, \emph{Arithmetic intersection theory},
Publ. Math. IHES \textbf{72}(1990), 94-174.

\bibitem{GS5}
H. Gillet and C. Soul\'{e}, \emph{An arithmetic Riemann-Roch
theorem}, Inventiones Math. \textbf{110}(1992), 473-543.

\bibitem{Hi}
F. Hirzebruch, \emph{Topological Methods in Algebraic Geometry},
Springer, 1978.

\bibitem{Hoe}
L. H\"{o}rmander, \emph{The analysis of linear partial
differential operators I}, Grundlehren der mathematischen
Wissenschaften \textbf{256}, Springer-Verlag, Berlin, 1977.

\bibitem{Ko}
B. K\"{o}ck, \emph{The Grothendieck-Riemann-Roch theorem for group
scheme actions}, Ann. Sci. Ecole Norm. Sup. \textbf{31}(1998),
4\`{e}me s\'{e}rie, 415-458.

\bibitem{KR1}
K. K\"{o}hler and D. Roessler, \emph{A fixed point formula of
Lefschetz type in Arakelove geometry I: statement and proof},
Inventiones Math. \textbf{145}(2001), no.2, 333-396.

\bibitem{Le}
P. Lelong, \emph{Int\'{e}gration sur un ensemble analytique
complexe}, Bull. Soc. Math. France \textbf{95}(1957), 239-262.

\bibitem{deRh}
G. deRham, \emph{Vari\'{e}t\'{e}s diff\'{e}rentiables: formes,
courants, formes harmoniques}, Hermann, Paris, 1966.

\bibitem{Th}
R. W. Thomason, \emph{Une formule de Lefschetz en K-th\'{e}orie
\'{e}quivariante alg\'{e}brique}, Duke Math. J. \textbf{68}(1992),
447-462.


\end{thebibliography}
\end{document}